\documentclass[11pt,citesort,rotate,psfig]{article}  
\usepackage{epsfig,amsfonts}
\catcode`\"=13      
\def"#1{\v{#1}}      

\def\be{\begin{equation}}
\def\ee{\end{equation}}
\def\bea{\begin{eqnarray}}
\def\eea{\end{eqnarray}}

\def\ppp#1 {#1^{\prime \prime \prime}}

\def\hfl#1#2{\smash{\mathop{\hbox to 12 mm{\rightarrowfill}} \limits^{\scriptstyle #1}_{\scriptstyle #2}}}

\def\binom#1#2{\left ( {{#1} \atop {#2}} \right )}

\newcommand {\qed} {\null \hfill \rule {2.5mm}{2.5mm}}
\newcommand {\R} {\ensuremath{\mathbb{R}}}
\newcommand {\N} {\ensuremath{\mathbb{N}}}

\author{Tomislav Do\v{s}li\'c$^{\star\dagger}$ and Darko Veljan$^{\star\ddagger}$ }
\title{Logarithmic behavior of some combinatorial sequences}
\date{\today}

\begin{document}
\maketitle
$\dagger$       Department of Informatics and Mathematics, Faculty of Agriculture, University of Zagreb, Sveto\v simunska c. 25, Zagreb, CROATIA
\\ \\
$\ddagger$        Department of Mathematics, University of Zagreb, Bijeni\v cka 30, Zagreb, CROATIA\\
\\
      $\star$ To whom correspondence should be addressed,
        e-mail : dveljan@cromath.math.hr\\

Proposed running head: Logarithmic behavior of combinatorial sequences

\newpage
{\bf \Large Abstract}
\\ \\
Two general methods for establishing the logarithmic behavior of recursively
defined sequences of real numbers are presented. One is the interlacing
method, and the other one is based on calculus. Both methods are used to
prove logarithmic behavior of some combinatorially relevant sequences, such as
Motzkin and Schr\"oder numbers, sequences of values of some classic
orthogonal polynomials, and many others. The calculus method extends also to
two- (or more- ) indexed sequences.

{\bf Keywords:} log-concavity, log-convexity, special combinatorial numbers, 
calculus, integer sequences

{\bf AMS subject classifications:} 05A20, 05A16, 05E35, 11B83, 11B37, 11B39,
11B68

\newpage
\section{Introduction}

Let $a(n)$, $n \geq 0$, be a sequence of positive real numbers. We 
want to examine the rate of growth of this sequence, i.e. to examine whether
the quotient $\frac {a(n)}{a(n-1)}$ decreases, increases or remains constant.
In other words, we want to see whether the sequence is {\bf log-concave}, i.e.
$a(n)^2 \geq a(n-1)a(n+1)$, {\bf log-convex}, i.e. $a(n)^2 \leq a(n-1)a(n+1)$,
or {\bf log-straight} (or {\bf geometric}), i.e. $a(n)^2=a(n-1)a(n+1)$ for 
all $n\geq 1$. Under log-behavior we also sometimes include {\bf log-Fibonacci}
behavior, meaning $sign[a(n)^2-a(n-1)a(n+1)]= sign (-1)^n$ (or $(-1)^{n+1}$).
It is of great interest, especially in combinatorics, as it
can be seen from many examples in \cite{stanley89}, to know the
log-behavior of a given sequence. It is, in fact, just one instance of the
whole paradigm of ``positivity questions'' (\cite{stanley00}). 

If $a(n)$ has a combinatorial meaning it would be ideally to provide a
combinatorial proof of its log-behavior. For example, if we want to prove that
$a(n)$ is log-convex and if we know that $a(n) = |S(n)|$, where $S(n)$ is a
certain finite set, then we would like to find an injection $S(n) \times
S(n) \rightarrow S(n-1) \times S(n+1)$, or a surjection $S(n-1) \times S(n+1)
\rightarrow S(n) \times S(n)$, and similarly for log-concavity. 
It is usually a hard task to find such a
(natural) injection or surjection. Still, examples of this type include
binomial coefficients, Motzkin numbers (\cite{callan}) and permutations with
a prescribed number of runs 
(\cite{bona}). Of course, the explicit formulae give another possibility to
prove results of this type, but they are rarely on disposal. Instead, other 
methods for proving such inequalities have been
developed, e.g. see \cite{sagan}, \cite{stanley89}, \cite{brenti} or \cite{bender}.

In this paper besides using old methods to prove some new results on
log-behavior, we shall also introduce some new methods and use them to prove
log-behavior of certain interesting combinatorial sequences, and apply this
method to other sequences, the most prominent example being values of
classical orthogonal polynomials.

\section{Log-behavior of some sequences using known results}

Let us quote some known results and apply them to examine the log-behavior of
certain combinatorial and other sequences.

{\bf Lemma 2.1} (Newton's lemma) \\ Let $P(x)=\sum_{k=0}^n a_k x^k $ be a real
polynomial whose all roots are real numbers. Then the coefficients of $P(x)$
are log-concave, i.e. $a_k^2 \geq a_{k-1}a_{k+1}$, $k=1, \ldots , n-1$.
Moreover, the (finite) sequence $\frac {a_k}{\binom {n}{k}}$ is log-concave
in $k$. \qed

Let us briefly recall how to apply this lemma to binomial coefficients and
Stirling numbers $c(n,k)$ of the first kind (the number of permutations on
the set $[n]=\{1, 2, \ldots , n\}$ with exactly $k$ cycles) and Stirling
numbers $S(n,k)$ of the second kind (the number of partitions of $[n]$ into
exactly $k$ blocks). The following formulae are well known:
\be
(x+1)^n=\sum_{k=0}^n\binom {n}{k} x^k,
\ee
\be
x^{\bar n} = \sum_{k=0}^n c(n,k) x^k,
\ee
\be
x^n=\sum_{k=0}^n S(n,k) x^{\underline {n}},
\ee
where $x^{\underline k}=x(x-1)\ldots (x-k+1)$ is the $k$-th falling power, and
$x^{\overline {k}}=x(x+1)\ldots (x+k-1)$ the $k$-th rising power of $x$. From
(2.1) and (2.2) we see that $(x+1)^n$ and $x^{\bar n}$ have only real roots. So,
by Newton's lemma we conclude that the sequences $\binom {n}{k}$ and $c(n,k)$
are log-concave. The case of the sequence $S(n,k)$ is a bit more involved. 
Let
$$P_n(x)=\sum_{k=0}^n S(n,k) x^k.$$
From $P_0(x)=1$ and from the basic recursion
$$ S(n,k)=S(n-1,k-1)+kS(n-1,k),$$
it follows at once that
$$P_n(x)=x[P_n'(x)+P_{n-1}(x)].$$
The function $Q_n(x)=P_n(x) e^x$ has the same roots as $P_n(x)$ and it is 
easy to verify $Q_n(x)=xQ_n'(x)$. By induction on $n$ and using Rolle's theorem
it follows easily that $Q_n$, and hence $P_n$, have only real and non-positive
roots. So, we conclude:

{\bf Theorem 2.2}\\ The sequences $\binom {n}{k} _{k \geq 0}$, $\left ( c(n,k)
\right ) _{k \geq 0}$ and $\left ( S(n,k)\right ) _{k \geq 0}$ are log-concave.
Hence, they are unimodal. \qed

An inductive proof of Theorem 2.2 is given in \cite{sagan}.

The next easy lemma is sometimes useful in proving log-convexity results.

{\bf Lemma 2.3}\\ Let $f:[a,b] \rightarrow \R$ be a positive, continuous function,
and $$I_n=\int _a ^b f(x)^n dx, \quad n \geq 1.$$ Then $\left ( I_n \right )
_{n \geq 2}$ is a log-convex sequence.

{\bf Proof}\\ By Cauchy-Schwarz inequality, we have
$$I_n^2=\left (\int _a ^b f(x)^n dx \right )^2 = 
\left (\int _a ^b f(x)^{\frac {n-1}{2}}f(x)^{\frac {n+1}{2}}dx \right )^2 
\leq \int _a ^b f(x)^{n-1}dx \int _a ^b f(x)^{n+1}dx 
=I_{n-1} I_{n+1}.
$$ \qed

As an example, we apply this lemma to Legendre polynomials $P_n(x)$. It
is well known (e.g. \cite{szego}) that the following Laplace formula holds:
\be
P_n(x)=\frac{1}{\pi }\int _0 ^{\pi} (x+\sqrt{x^2-1} cos \varphi )^n d 
\varphi.
\ee
Hence, from Lemma 2.3 and (2.4) we obtain

{\bf Theorem 2.4}\\ The values $P_n(x)$, $n \geq 0$, for $x \geq 1$ are 
log-convex. \qed

Another proof of this fact will be presented in Section 4.

We say that a sequence $\left ( a_n \right ) _{n \geq 0}$ has no internal zeros
if there do not exist integers $0 \leq i < j<k$ such that $a_i\neq0$, $a_j=0$,
$a_k\neq 0$.

{\bf Theorem 2.5} (Bender-Canfield, see \cite{bender}) \\ Let $1, a_1, a_2,
\ldots $ be a log-concave sequence of nonnegative real numbers with no 
internal zeros and let $\left
( b_n \right ) _{n \geq 0}$ be the sequence defined by
\be
\sum _{n \geq 0} b_n \frac{x^n}{n!} = exp \left ( \sum _{k \geq 1} a_k
\frac{x^k}{k!} \right ).
\ee
Then the sequence $(b_n)_{n\geq 0}$ is log-convex and $\left ( \frac{b_n}{n!}
\right ) _{n\geq 0}$ is log-concave. \qed

As our first application of this theorem consider the Bell numbers $(B_n)_{n
\geq 0}$. $B_0=1$, and $B_n$ is the number of partitions of an $n$-set. It is
well known that the exponential generating function for $\left ( B_n \right )_{
n \geq 0}$ is given by
$$\sum _{n \geq 0} B_n \frac {x^n}{n!} = exp (e^x-1) = exp \left (
\sum _{k \geq 1} \frac {x^k}{k!} \right ) .$$
By taking $a_k = \frac{1}{(k-1)!}$, $k \geq 1$, and checking that the sequence
$1, a_1, a_2, \ldots $ is log-concave, we conclude from Theorem 2.5 that
$\left ( B_n \right )_{n \geq 0}$ is log-convex and $\left ( \frac{B_n}{n!}
\right )_{n \geq 0}$ is log-concave sequence. More generally, for an integer
$l \geq 2$, define $exp_l$ to be $l$ times iterated exponential function, i.e.
$$exp_l(x)=exp(exp(\ldots (exp(x))\ldots ),$$
$exp$ written $l$ times. Define the sequence $\left ( b_n^{(l)} \right )_{n 
\geq 0}$ by
$$\sum _{n \geq 0} b_n^{(l)} \frac {x^n}{n!}= exp_l(x),$$
and {\bf Bell numbers of order $l$} by
\be
B_n^{(l)} = \frac{b_n^{(l)}}{exp_l(0)}.
\ee
So, $B_n^{(2)} = B_n$ are ordinary Bell numbers. Now it is not hard to prove
by induction on $l$ the following result.

{\bf Theorem 2.6}\\ For any fixed $l \geq 2$, the Bell numbers of order $l$, i.e.
the sequence $\left ( B_n^{(l)} \right )_{n\geq 0}$ is log-convex, and the
sequence $\left ( \frac{B_n^{(l)}}{n!}\right )_{n\geq 0}$ is log-concave. \qed

The following lemma is an easy consequence of the definition of log-convexity
and log-concavity.

{\bf Lemma 2.7}\\ Suppose $\left ( a_n \right )_{n\geq 0}$ is a positive, 
log-convex sequence and $a_0 = 1$. Then $a_na_m \leq a_{n+m}$. If, in addition,
$\left ( \frac{a_n}{n!}\right )_{n \geq 0}$ is log-concave, then $a_{n+m}
\leq \binom{m+n}{n}a_na_m$. \qed

Now Theorem 2.6 and Lemma 2.7 imply semi-additivity inequalities for Bell numbers
of order $l$.

{\bf Corollary 2.8} $$B_n^{(l)} B_m^{(l)} \leq B_{m+n}^{(l)} \leq \binom{m+n}{n}
B_n^{(l)} B_m^{(l)}, \quad m,n \geq 0.$$ \qed

The next application of the Bender-Canfield theorem concerns the number of
certain permutations, i.e. elements $\pi \in \Sigma _n$, the symmetric group of
permutations on $n$ letters. For a fixed integer $k \geq 1$ let $c_k(n)$ be
the number of all permutations from $S_n$ that have cycles of length not
exceeding $k$. By definition, $c_k(0)=1$.

{\bf Theorem 2.9}\\ The sequence $\left ( c_k(n)\right )_{n\geq 0}$ is log-convex
and the sequence $\left ( \frac{c_k(n)}{n!}\right )_{n \geq 0}$ is log-concave.

{\bf Proof}\\ It is well known (see, e.g. \cite{stanleyII}) that the
exponential generating function $ C_k(x)$ of the sequence $\left ( 
c_k(n) \right )_{n\geq 0}$ is given by 
$$C_k(x)= exp \left ( \sum _{j=1}^k \frac{x^j}{j} \right ). $$
Hence the sequence $\left ( a_k\right )_{k\geq 0}$ from the Bender-Canfield
theorem in (2.5) is as follows: $a_0=1, a_1=1, \ldots , a_k=1, a_{k+1}=0,
a_{k+2}=0, \ldots $, which has no internal zeros and is obviously a log-concave sequence. Hence the 
claim follows from Theorem 2.6. \qed

The requirement that the sequence $(a_n)$ from Bender-Canfield theorem does not
have internal zeros is essential, and in general can not be weakened. As an
illustration, let us consider a class of sequences related to
$c_k(n)$.
Again, for a fixed integer $k \geq 1$ let $e_k(n)$ be the number of all
permutations $\pi $ from $\Sigma _n$ such that $\pi ^k =id$, for $n \geq 1$ and
$e_k(0)=1$. Then the exponential generating function $E_k(x)$ of this sequence
is given by (see, e.g. \cite{stanleyII})
\be
E_k(x)= exp \left ( \sum _{j |k} \frac{x^j}{j} \right ), 
\ee
and the corresponding recurrence is 
$$
e_k(n)=\sum _{j | k} (n-1)^{\underline {j-1}} e_k(n-j),
$$
with appropriate initial conditions. The sequence $(a_n)$ from (2.5) is given
by $a_0=a_1=1$, $a_j=1$ if $j$ divides $k$ and $a_j=0$ otherwise. This binary
sequence $1,1,a_2, \ldots , a_k,0,0,\ldots $ is log-concave if and only if it
does not contain $1,0,1$ as a subsequence, but it contains internal zeros for
all $k > 2$. Taking, for example, $k=5$, the first few terms of $e_5(n)$
being $1,1,1,1,1,25,145,505,1345,\ldots $, we can easily see that this sequence
does not exhibit any logarithmically definite behavior, although the sequence
$1,1,0,0,0,1,0,\ldots $ is log-concave. For higher values of $k$, sequences
$e_k(n)$ log-behave even more chaotically.

\section{The interlacing (or ``sandwich'') method. Secondary structures.}

\setcounter{equation}{0}
Let $a(n)$, $n \geq 0$, be a sequence of positive numbers defined by a
homogeneous linear recurrence, say of second order:
\be
a(n)=R(n)a(n-1)+S(n)a(n-2),
\ee
where $R$ and $S$ are known functions, together with given initial values
$a(0)=a_0$, and $a(1)=a_1$.

Let our task be to examine the rate of growth of $a(n)$. We define the sequence
of consecutive quotients
$$q(n)=\frac{a(n)}{a(n-1)}, \quad n \geq 1.$$
Dividing (3.1) by $a(n-1)$ we obtain the recurrence
\be
q(n) = R(n) + \frac{S(n)}{q(n-1)},
\ee
with initial condition $q(1)=\frac{a_1}{a_0}:=b_1$. The log-concavity or 
log-convexity of $(a(n))$ is equivalent, respectively, to $q(n) \geq q(n+1)$
or $q(n) \leq q(n+1)$, for all $n \geq 1$. So, what we want to see is whether 
the sequence $\left ( q(n) \right ) _{n \geq 1}$ decreases or increases. To
prove that $(q(n))$ increases, it is enough to find an increasing sequence
$(b(n))$ such that 
\be
b(n) \leq q(n) \leq b(n+1)
\ee
holds for all $n \geq 1$, or at least for all $n \geq n_0$ for some $n_0$.
Then we can conclude that $(a(n))$ behaves log-convex at least from some
place on. Analogously for log-concavity. This ``sandwich method'' or 
``interlacing method'' works in some simple cases, but often it is very hard
to hit the right sequence $(b(n))$ which is simple enough. In the rest of this
section we show how this method works for some combinatorially important
sequences. We also show some consequences of the obtained results.

{\bf Example 3.1} (Derangements) \\ Let $D_n$ be the number of derangements
on $n$ objects, i.e. the number of permutations $\pi \in \Sigma _n$ without fixed
points, for $n \geq 1$, and $D_0 := 1$. It is well known (and easy to prove)
that the following recurrence holds:
\be
D_n = (n-1)[D_{n-1}+D_{n-2}],
\ee
with initial conditions $D_0=1$, $D_1=0$. Then $D_2=1$, $D_3=2$, $D_4=9$, 
$D_5=44$, $D_6=265$ etc., and we expect from these initial values that
$D_n^2 \leq D_{n-1}D_{n+1}$ for $n \geq 3$. Indeed, divide (3.4) by $D_{n-1}$
and denote $q(n)=\frac{D_n}{D_{n-1}}$. Then
\be 
q(n)=(n-1)\left [ 1+ \frac{1}{q(n-1)} \right ],
\ee
with $q(3)=2$, $q(4)=9/2$. Let $b(n)=n-1/2$. It is easy to check by induction
on $n$ and using (3.5) that
$$ b(n) \leq q(n) \leq b(n+1), $$
for all $n \geq 4$. Since $(b(n))$ is clearly increasing and $q(3) \leq q(4)$,
we conclude  that $(D_n)$ is log-convex for $n \geq 3$. \qed

\newpage
{\bf Example 3.2} \\
Let $T_2(n)$ denotes the number of $n \times n$ symmetric $\N _0$-matrices with
every row (and hence every column) sum equal to $2$ with trace zero (i.e. all
main-diagonal entries are zero) (Example 5.2.8 in \cite{stanleyII}). The
exponential generating function of $T_2(n)$ is given by
$$T(x)=\frac{1}{\sqrt {1-x}}exp\left ( \frac{x^2}{4}-\frac{x}{2} \right ) .$$
The numbers $T_2(n)$ satisfy the recurrence
$$T_2(n)=(n-1)T_2(n-1)+(n-1)T_2(n-2)-\binom{n-1}{2}T_2(n-3)$$
with the initial conditions $T_2(0)=1$, $T_2(1)=0$, $T_2(2)=T_2(3)=1$.
The corresponding recurrence for successive quotients $q_2(n)=T_2(n)/T_2(n-1)$
is given by
$$q_2(n)=(n-1)+\frac{n-1}{q_2(n-1)}-\binom{n-1}{2}\frac{1}{q_2(n-1)q_2(n-2)},
\quad n \geq 5,$$
with the initial conditions $q_2(3)=1$, $q_2(4)=6$. Tabulating the first few
values of $q_2(n)$, we see that, after some initial fluctuations, this 
sequence seems to behave like $n-\frac{1}{2}$. So, we guess that $n-1 \leq
q_2(n) \leq n$, and indeed, this follows easily by induction on $n$ for $n
\geq 6$. Hence, the sequence $\left ( T_2(n) \right ) _{n \geq 6}$ is 
log-convex. \qed

For our next application we need some preparations. A {\bf Motzkin path} of length
$n$ is a lattice path in $(x,y)$-plane from $(0,0)$ to $(n,0)$ with steps
$(1,1)$ (or $Up$), $(1,-1)$ (or $Down$) and $(1,0)$ (or $Level$), never 
falling below the $x$-axis. Denote by ${\cal M}(n)$ the set of all Motzkin paths
of length $n$. The number $M_n=|{\cal M}(n)|$ is the $n$-th {\bf Motzkin number}.
By definition, $M_0=1$.

A handful of other combinatorial interpretations of $M_n$ are listed in Ex.
6.38 in \cite{stanleyII}. A typical member of the Motzkin family ${\cal M}(20)$
is shown in Fig. 1 as a ``landscape path''.
\begin{figure}[h] \centerline {
\epsfig{file=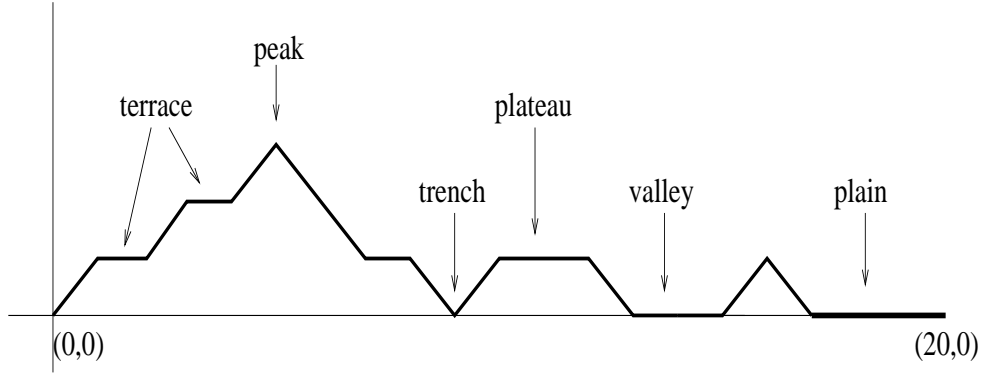,height=5cm,width=13cm,silent=}} 
\caption{A Motzkin path of length 20}
\label{motz}
\end{figure}
A {\bf peak} of a Motzkin path is a place where an $Up$ step is immediately
followed by a $Down$ step. A {\bf plateau} of length $l$ is a sequence of $l$
consecutive $Level$ steps, immediately preceded by an $Up$ step and
immediately followed by a $Down$ step. Similarly we define a {\bf terrace},
{\bf trench}, {\bf valley} and {\bf plain} of a Motzkin path.

A Motzkin path without any $Level$ steps is called a {\bf Dyck path} or a ``mountain path''. It is 
well known that the set ${\cal D}(n)$ of all Dyck paths of length $2n$ is
enumerated by {\bf Catalan numbers} $C_n$, i.e. $|{\cal D}(n)| = C_n =
\frac{1}{n+1}\binom{2n}{n}$. From the explicit formula
it follows immediately that $C_n^2 \leq C_{n-1}C_{n+1}$, i.e. Catalan numbers
are log-convex. It is also easy to find a simple combinatorial proof of this
fact. By counting the number of $Level$ steps on a Motzkin path, it can be
easily shown that Catalan and Motzkin numbers are related as follows:
$$M_n=\sum _{k \geq 0}\binom{n}{2k}C_k, \quad
C_{n+1}=\sum _{k \geq 0}\binom{n}{k}M_k.$$

{\bf Proposition 3.2}\\
a) The Motzkin numbers satisfy the following convolutive recursion:
\be
M_{n+1}=M_n+\sum _{k =0}^{n-1} M_k M_{n-k-1}.
\ee
b) The generating function of $(M_n)$ is given by
\be
M(x)=\sum _{n \geq 0} M_n x^n =\frac{1-x-\sqrt{1-2x-3x^2}}{2x^2}.
\ee
c) $M_n$'s satisfy the short recursion
\be
(n+2)M_n=(2n+1)M_{n-1}+3(n-1)M_{n-2}.
\ee

{\bf Proof}\\ a) $M_0=1$ and $M_1=1$. A Motzkin path of length $n+1$ either
starts by a $Level$ step and then can proceed in $M_n$ ways, or starts by an
$Up$ step and returns for the first time to the $x$-axis after $k+1$ steps
(to the point $(k+2,0)$). The number of latter is equal to the number of pairs
$(P_1,P_2)$, where $P_1$ is a (translated) Motzkin path from $(1,1)$ to $(
k+1,1)$ which is not below the line $y=1$, and $P_2$ is a (translated) Motzkin
path from $(k+2,0)$ to $(n+1,0)$. The number of paths $P_1$ is equal to $M_k$,
while the number of paths $P_2$ is equal to the number of Motzkin paths on
$n+1-(k+2)=n-k-1$ steps, and this is $M_{n-k-1}$. Thus a) follows.

b) If we multiply (3.6) by $x^{n+1}$ and sum over $n \geq 0$, we get the
functional equation $x^2 M^2(x)+(x-1)M(x)+1=0$, and since $M(0)=M_0=1$, we
obtain (3.7).

c) The generating function $M(x)$ is algebraic, hence D-finite. 
Therefore $(M_n)$ is $P$-recursive. Now from Eq(6.38) in
\cite{stanleyII} for the polynomial $A(x)=1-2x-3x^2$ and $r=2$, $a_0=1$, 
$a_1=-2$, $a_2=-3$, $d=2$, we get the claim. \qed

{\bf Theorem 3.3} \\
a) The sequence $M_n$ of Motzkin numbers is log-convex.\\
b) $M_n/M_{n-1} < 3$, for all $n \geq 1$. \\
c) The sequence $M_n/M_{n-1}$ is convergent and $\lim _{n \rightarrow \infty}
M_n/M_{n-1} = 3$.

\newpage
{\bf Proof}\\
a) Divide (3.8) by $(n+2)M_{n-1}$ and let $q(n):=M_n/M_{n-1}$. Then we obtain
\be
q(n)=\frac{1}{n+2}\left [ 2n+1 + \frac{3(n-1)}{q(n-1)}\right ], \quad n\geq 2,
\ee
with the initial condition $q(1)=1$. Then $q(2)=2$. Now define the sequence
$b_n = 3\frac{2n}{2n+3}$. This sequence is obviously increasing. We claim that
$b_n \leq q(n) \leq b_{n+1}$, for all $n \geq 3$. This is obviously true for
$n=3$. Let $n \geq 3$. By inductive hypothesis, we have
$$q(n)=\frac{2n+1}{n+2}+\frac{3(n-1)}{n+2}\frac{1}{q(n-1)} \geq
\frac{2n+1}{n+2}+\frac{3(n-1)}{n+2}\frac{1}{b_n}.$$
But
$$\frac{2n+1}{n+2}+\frac{3(n-1)}{n+2}\frac{1}{b_n}-b_n = \frac{3}{2} 
\frac{n-3}{n(n+2)(2n+3)} \geq 0,$$
for all $n \geq 3$, and the inequality $q(n) \geq b_n$ follows. On the other 
hand, by the inductive hypothesis $b_{n-1} \leq q(n-1)$, we have
$$q(n) \leq \frac{2n+1}{n+2}+\frac{3(n-1)}{n+2}\frac{1}{b_{n-1}}.$$
Subtracting $b_{n+1}$ from the right hand side, we get
$$\frac{2n+1}{n+2}+\frac{3(n-1)}{n+2}\frac{1}{b_{n-1}}-b_{n+1} = 
-\frac{1}{2}\frac{1}{(n+2)(2n+5)} \leq 0,$$
and the inequality $q(n) \leq b_{n+1}$ follows. So, the claim is proved by
induction. Hence, the sequence $M_n/M_{n-1}$ is increasing, and the Motzkin
sequence is log-convex.

b) We have proved in a) that $M_n/M_{n-1} \leq b_n = 3\frac{2n}{2n+3} < 3$.

c) By a) and b) it follows that $q(n)=M_n/M_{n-1}$ is an increasing and
bounded sequence: $2 \leq q(n) < 3$, hence convergent. Passing to limit in
(3.9), or passing to limit in the sandwich inequality $b_n \leq q(n) \leq
b_{n+1}$ above, the last claim in c) follows. \qed

{\bf Corollary 3.4}\\
a) The sequence $(M_n/n!)$ is log-concave. \\
b) $M_n^2 \leq M_{n-1}M_{n+1} \leq \left ( 1+\frac{1}{n} \right )M_n^2$, for 
all $n \geq 1$.\\
c) $M_mM_n \leq M_{m+n} \leq \binom{m+n}{n}M_mM_n$, for all $m,n \geq 0$.

{\bf Proof}\\
a) The log-concavity of $M_n/n!$ is equivalent to $\frac{M_{n+1}}{M_n} \leq
\frac{n+1}{n}\frac{M_n}{M_{n-1}}$. So, we have to prove that $q(n+1) \leq
\frac{n+1}{n} q(n)$. We know that the sequence $q(n)$ is increasing. Starting
from the short recursion (3.9) for $q(n+1)$, we have
\begin{eqnarray*}
q(n+1)& =& \frac{2n+3}{n+3}+\frac{3n}{n+3}\frac{1}{q(n)} \leq
\frac{2n+3}{n+3}+\frac{3n}{n+3}\frac{1}{q(n-1)} \\
& =&
\frac{n+2}{n+3}\frac{n}{n-1}\left [ \frac{(2n+3)(n-1)}{2n(n+1)}\frac{2n+1}{n+2}
+\frac{3(n-1)}{n+2}\frac{1}{q(n-1)}\right ].
\end{eqnarray*}
The claim now follows by noting that the term $\frac{(2n+3)(n-1)}{2n(n+1)}$
is clearly less than one for all $n \geq 1$, and that the inequality
$\frac{n+2}{n+3}\frac{n}{n-1} < \frac{n+1}{n}$ is valid for all $n \geq 3$.
The validity of our claim can easily be checked for $1 \leq n \leq 3$.

b) This follows from log-convexity of Motzkin numbers and a). 

c) A simple combinatorial proof of the left inequality follows from the fact 
that a concatenation of two Motzkin paths of lengths $m$ and $n$, respectively,
is again a valid Motzkin path of length $m+n$. To prove the right inequality,
start from $q(n) \geq \frac{n}{n+1} q(n+1)$. Using this inequality repeatedly,
we get
$$\frac{M_1}{M_0} \geq \frac{1}{2}\frac{M_2}{M_1} \geq \frac{1}{3}\frac{M_3}
{M_2} \geq \ldots \geq \frac{1}{m+n}\frac{M_{m+n}}{M_{m+n-1}},$$
for all $n \geq 0$, $m \geq 1$.

Hence, for any $0 \leq j \leq m-1$, we have
$$\frac{M_{j+1}}{M_j}\geq \frac{j+1}{m+n}\frac{M_{m+n}}{M_{m+n-1}}.$$
From this we get
$$\frac{M_1}{M_0}\frac{M_2}{M_1} \ldots \frac{M_m}{M_{m-1}} \geq
\left ( \frac{1}{n+1}\frac{M_{n+1}}{M_n}\right ) 
\left ( \frac{2}{n+2}\frac{M_{n+2}}{M_{n+1}} \right ) \ldots 
\left ( \frac{m}{m+n}\frac{M_{m+n}}{M_{m+n-1}} \right ),$$
and after the cancellation:
$$\frac{M_m}{M_0} \geq \frac{m!n!}{(m+n)!}\frac{M_{m+n}}{M_n}.$$
Since $M_0=1$, we get the claim. The case $m=0$, $n \geq 0$ is trivial. \qed

{\bf Remark 3.5}\\ An algebraic proof of log-convexity of Motzkin numbers was
given in \cite{aigner}, and a combinatorial proof in \cite{callan}. We shall
give yet another (``calculus'') proof later in Section 4.

Now we proceed in applying the ``sandwich method'' to combinatorial structures
that generalize both Dyck and Motzkin structures in a sense that they are
counted by Catalan and Motzkin numbers in some special cases. They are called
secondary structures and come from molecular biology (see, e.g. \cite{stwat},
\cite{waterm}, \cite{watermhand}, \cite{stadler}, \cite{hofack}, \cite{kruskal}). More details are in 
\cite{dsv} and \cite{doslicphd}.

Let $n$ and $l$ are integers, $n\geq 1$, $l \geq 0$. A {\bf secondary structure}
of {\bf size} $n$ and {\bf rank} $l$ is a labeled graph $S$ on the vertex set
$V(S)=[n]=\{1,2,\ldots ,n\}$ whose edge set $E(S)$ consists of two disjoint
subsets, $P(S)$ and $H(S)$, satisfying the following conditions:\\
(a) $\{i,i+1\} \in P(S)$, for all $1 \leq i \leq n-1$;\\
(b) $\{i,j\} \in H(S)$ and $\{i,k\} \in H(S) \Longrightarrow j=k$;\\
(c) $\{i,j\} \in H(S)\Longrightarrow |i-j|>l$;\\
(d) $\{i,j\} \in H(S), \{p,q\} \in H(S)$ and $i<p<j$ $\Longrightarrow i<q<j$.

The vertices of a secondary structure are (in biology) usually called {\bf
bases}, the edges in $P(S)$ are called {\bf $p$-bonds} and those in $H(S)$
{\bf $h$-bonds}. A secondary structure $S$ with $H(S)=\emptyset$ is called
{\bf trivial}. The number of $h$-bonds of $S$ is called its {\bf order}.
An example of a secondary structure of size 12, rank 2 and order 3 is shown
below in Fig. 2. Note that every $h$-bond ``leaps'' over at least two bases.
\begin{figure}[h] \centerline {
\epsfig{file=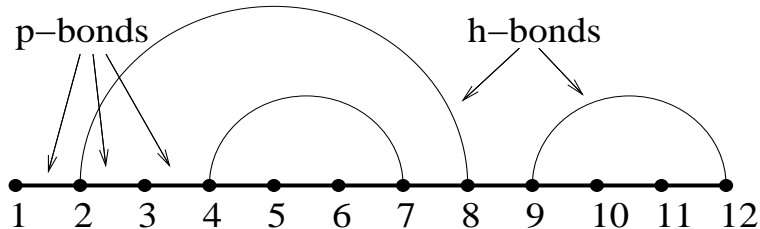,height=3cm,width=10cm,silent=}}
\caption{An example of a secondary structure}
\label{sec}
\end{figure}
We denote the set of all secondary structures of size $n$ and rank $l$ by
${\cal S}^{(l)}(n)$, and the set of such structures of order $k$ by$ {\cal S}
^{(l)}_k(n)$. The respective cardinalities we denote by $S^{(l)}(n)$ and
$S^{(l)}_k(n)$. By definition, we put $S^{(l)}(0)=1$, for all $l$.

Another interpretation of secondary structures is as follows. Denote by
${\cal M}^{(l)}(n)$ the set of all Motzkin paths in ${\cal M}(n)$ whose every
plateau is at least $l$ steps long. For $l=0$, ${\cal M}^{(0)}(n)={\cal M}(n)$,
i.e. all Motzkin paths on $n$ steps.

{\bf Proposition 3.6}\\ There is a bijection between ${\cal S}^{(l)}(n)$ and
${\cal M}^{(l)}(n)$ for all $n\geq 1$, $l \geq 0$.

{\bf Proof}\\ We shall only briefly describe the correspondence ${\cal S}^{(l)}
(n) \rightarrow {\cal M}^{(l)}(n)$. Take a secondary structure $S \in
{\cal S}^{(l)}(n)$ and scan it from left to right. To each base, starting
from $1$, we assign a step in a lattice path as follows. To any unpaired base
(i.e. base in which no $h$-bond starts or ends) we assign a $Level$ step. To a
base in which a new $h$-bond starts assign an $Up$ step, and to each base in
which an already encountered $h$-bond terminates, assign a $Down$ step. The
obtained path $P$ is in ${\cal M}^{(l)}(n)$ and it is not hard to check that
$S \longmapsto P$ is a bijection. \qed

An example of the correspondence ${\cal S}^{(1)}(10) \longleftrightarrow
{\cal M}^{(1)}(10)$ is shown in Fig. 3.
\begin{figure}[h] \centerline {
\epsfig{file=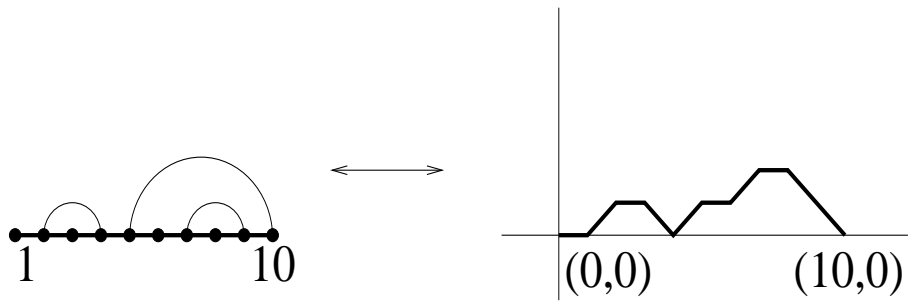,height=4cm,width=12cm,silent=}} 
\caption{A secondary structure and the corresponding Motzkin path}
\label{corr}
\end{figure}
Therefore, the rank $0$ secondary structures, i.e the border case, are just
all Motzkin paths, hence $S^{(0)}(n)=M_n$. Note that for $l \geq 1$ a secondary
structure of rank $l$ is a simple graph.

If we allow rank to degenerate to the value of $l=-1$, then an $h$-bond can
terminate in the very vertex it starts from, i.e. we allow loops.

{\bf Proposition 3.7}\\ There is a bijection between the set ${\cal S}^{(-1)}
(n)$ of all secondary structures of size $n$ and rank $-1$ and the set 
${\cal D}(n+1)$ of all Dyck paths on $2(n+1)$ steps. Hence, $S^{(-1)}(n)=
C_{n+1}$, for all $n \geq 0$.

{\bf Proof}\\ Again, we shall only briefly describe how to assign a member of
${\cal S}^{(-1)}(n)$ to a Dyck path $P \in {\cal D}(n+1)$. So, take a Dyck
path on $2(n+1)$ steps. Discard the first and the last step and divide the
remaining path in groups of two consecutive steps. Assign to each group a base
in a secondary structure according to the following rule. To a group of two
$Up$ steps assign a base in which an $h$-bond starts, to a group of two $Down$
steps assign a base in which the last started $h$-bond terminates. To a group
$(Up,Down)$ assign a base with a loop attached to it, and, finally, to a group
$(Down,Up)$ assign an unpaired base. The obtained graph $S$ is in 
${\cal S}^{(-1)}(n)$ and $P \longmapsto S$ is a bijection. \qed

An example of the correspondence ${\cal D}(8) \longleftrightarrow {\cal S}
^{(-1)}(7)$ is shown in Fig. 4.
\begin{figure}[h] \centerline {
\epsfig{file=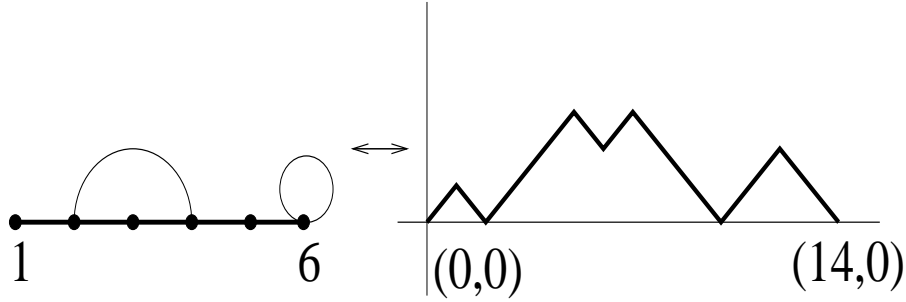,height=4cm,width=12cm,silent=}} 
\caption{The correspondence between secondary structures of rank -1 and Dyck paths}
\label{lmin1}
\end{figure}
The third combinatorial interpretation of secondary structures is related to
pattern avoidance in permutations. A {\bf pattern} is a permutation $\sigma \in
\Sigma _k$, and a permutation $\pi \in \Sigma _n$ {\bf avoids} $\sigma $ if
there is no subsequence in $\pi $ whose members are in the same relative order
as the members of $\sigma $. It is well known (\cite{knuth}) that the number
of permutations from $\Sigma _n$ avoiding $\sigma \in \Sigma _3$ is equal to
the Catalan number $C_n$, for all patterns $\sigma \in \Sigma _3$.

The concept of pattern avoidance was generalized in \cite{babson}, by allowing
the requirement that two letters adjacent in a pattern must be adjacent in
the permutation. An example of a generalized pattern is $1-32$, where an
$1-32$ sub-word of a permutation $\pi = a_1a_2\ldots a_n$ is any sub-word
$a_ia_ja_{j+1}$ such that $i<j$ and $a_i < a_{j+1} < a_j$. Generalized pattern
avoidance is treated in more detail in \cite{claesson}, where it is shown that
the permutations from $\Sigma _n$ that avoid both $1-23$ and $13-2$ are
enumerated by the Motzkin numbers.

{\bf Proposition 3.8}\\ For any $l \geq 0$, there is a bijection between ${\cal S}^{(l)}(n)$ and 
the set of all permutations from $\Sigma _n$ that avoid $\{1-23,13-2,ij\}$,
where $j \geq i+l$.

{\bf Proof}\\ This bijection is just the bijection given in the proof of
Proposition 24 in \cite{claesson}. It is easy to see that any pattern $i,i+k$
in a permutation $\sigma \in \Sigma _n$ 
avoiding $\{1-23,13-2\}$ generates a plateau of length $k-1$
in the corresponding Motzkin path. So, the claim follows from Proposition
3.6. \qed

Let us find now the recurrences and generating function for the general
secondary structure numbers analogous to those for the Motzkin numbers (the
border case) as was shown in Proposition 3.2.

{\bf Proposition 3.9}\\ 
a) For any fixed integer $l \geq -1$, the numbers $S^{(l)}(n)$ of secondary
structures of rank $l$ satisfy the following convolutive recurrence
\be
S^{(l)}(n+1)=S^{(l)}(n)+\sum _{m=l}^{n-1}S^{(l)}(m)S^{(l)}(n-m-1), \quad
n \geq l+1,
\ee
together with the initial conditions
\be
S^{(l)}(0)=S^{(l)}(1)=\ldots =S^{(l)}(l+1)=1.
\ee
b) The generating function $S^{(l)}(x)=\sum _{n\geq 0} S^{(l)}(n) x^n$ is given
by
\be
S_l(x)=\frac{-\omega _l(x)-\sqrt{\omega _l^2(x)-4x^2}}{2x^2},
\ee
where
\be
\omega _l(x)=x-1-x^2- \ldots -x^{l+1} = 2x-(1+x+x^2+ \ldots + x^{l+1}) = 
1-x^2\frac{1-x^l}{1-x}.
\ee
c) The sequence $\left ( S^{(l)}(n) \right )_{n\geq 0}$ satisfies the following
``short'' recursion
\be
(n+2)S^{(l)}(n)=\sum _{k=1}^{2l+2}A^{(l)}(n,k)S^{(l)}(n-k),
\ee
where
$$A^{(l)} (n,k)= \left\{\matrix{
\hfill -\frac{1}{2} (k-3)(2n+4-3k) &,& \hfill 1 \leq k \leq l+1 \cr
\hfill -\frac{1}{2} (l-3)(2n-3l-2) &,& \hfill k=l+2 \cr
-\frac{1}{2} (2l+3-k)(2n+4-3k) &,& l+3 \leq k \leq 2l+2 .\cr
} \right. $$

{\bf Proof}\\
a) Clearly, for $n \leq l+1$ there is only one (the trivial) secondary structure
of rank $l$ and size $n$, and hence initial conditions (3.11) hold. Let $n \geq
l+1$. A secondary structure on $n+1$ bases either does not contain an $h$-bond
starting at the base $1$, in which case there are $S^{(l)}(n)$ such structures,
or it has an $h$-bond from the base $1$ to some base $m+2$, at least $l$ bases
apart. In this case, there are $S^{(l)}(m)S^{(l)}(n-m-1)$ such structures, and
(3.10) follows.

b) Multiplying (3.10) by $x^{n+1}$, summing over $n \geq 0$ and taking into
account the initial conditions (3.11), we obtain the functional equation
$$x^2[S_l(x)]^2+\omega _l(x) S_l(x)+1=0,$$
and this, in turn, implies (3.12).

c) As in the proof of Proposition 3.2 c), we use D-finiteness of the above
generating function and again formula (6.38) in \cite{stanleyII} with $d=2$, $r=2l+2$
to obtain the claim. More details are in \cite{doslicphd} and \cite{dsv}. \qed

Let us write down explicitly the recurrences (3.14) in cases $l=1,2,3$ (of
course, for $l=0$, it coincides with (3.8)):
\be
S^{(1)}(n)=\frac{2n+1}{n+2}S^{(1)}(n-1)+\frac {n-1}{n+2}S^{(1)}(n-2)+
\frac {2n-5}{n+2}S^{(1)}(n-3)- \frac{n-4}{n+2}S^{(1)}(n-4),
\ee
$$S^{(1)}(0)=S^{(1)}(1)=S^{(1)}(2)=1, S^{(1)}(3)=2;$$

\be
S^{(2)}(n)=\frac{2n+1}{n+2}S^{(2)}(n-1)+\frac {n-1}{n+2}S^{(2)}(n-2)+
\frac{n-4}{n+2}S^{(2)}(n-4)- \\ \frac {2n-11}{n+2}S^{(2)}(n-5)- \frac{n-7}
{n+2}S^{(2)}(n-6),
\ee
$$
S^{(2)}(0)=S^{(2)}(1)=S^{(2)}(2)=S^{(2)}(3)=1, S^{(2)}(4)=2, S^{(2)}(5)=4;
$$

\bea
\lefteqn{S^{(3)}(n)=\frac{2n+1}{n+2}S^{(3)}(n-1)+\frac {n-1}{n+2}S^{(3)}(n-2)-
\frac{n-4}{n+2}S^{(3)}(n-4)} \nonumber \\
& & - \frac {3n-21}{n+2}S^{(3)}(n-6)-\frac {2n-17}
{n+2}S^{(3)}(n-7)-\frac {n-10}{n+2}S^{(3)}(n-8), 
\eea
$$S^{(3)}(0)=S^{(3)}(1)=S^{(3)}(2)=S^{(3)}(3)=S^{(3)}(4)=1,S^{(3)}(5)=2,
S^{(3)}(6)=4,S^{(3)}(7)=8.$$

{\bf Remark 3.10}\\ The recurrences (3.14) do not have any polynomial solutions
or solutions in hypergeometric terms for $l \geq 0$. This fact follows by
applying the algorithm {\bf Hyper}, described in \cite{petkovsek}.

Finally, we now return to our main theme: the log-convexity of secondary 
structure numbers.

{\bf Theorem 3.11}\\ The sequence $\left ( S^{(1)}(n) \right ) _{n \geq 0}$
is log-convex.

{\bf Proof}\\ We start from the short recursion (3.15). Dividing (3.15) by
$S^{(1)}(n-1)$ and denoting $S^{(1)}(n)/S^{(1)}(n-1)$ by $q_n$, we obtain
the following recursion for the numbers $q_n$:
\be
q_n=\frac{1}{n+2}\left [ 2n+1+\frac{n-1}{q_{n-1}}+\frac{2n-5}{q_{n-1}q_{n-2}}-
\frac{n-4}{q_{n-1}q_{n-2}q_{n-3}} \right ], \quad n\geq 4,
\ee
with the initial conditions $q_1=q_2=1$, $q_3=2$. It is easy to check that
$q_1 \leq q_2 \leq \ldots \leq q_6$.

Assume, for the moment (and we shall prove it later on) that the sequence 
$(q_n)$ is convergent with limit $q$ when $n \rightarrow \infty$. By passing to
limit in (3.18), we obtain the equation $q^4-2q^3-q^2-2q+1=0$, whose maximal
positive solution is $q=\varphi ^2=\varphi +1$, where $\varphi =\frac{1+\sqrt 5}
{2}$ is the golden ratio.

Define now the sequence $a_n=\frac{2n}{2n+3}\varphi ^2$. It is clearly an
increasing sequence and its limit is $\varphi ^2$. We claim that $(a_n)$ is
interlaced with our sequence $(q_n)$. More precisely, we shall prove by 
induction that
\be
a_n \leq q_n \leq a_{n+1}, 
\ee
for $n \geq 6$.

First we check directly the cases $n=6,7,8$ and $9$. Now take $n \geq 9$. From
the induction hypothesis and (3.18), we have
\begin{eqnarray*}
(n+2)q_n &=& 2n+1+\frac{n-1}{q_{n-1}}+\frac{n-1}{q_{n-1}q_{n-2}}+
\frac{n-4}{q_{n-1}q_{n-2}}-\frac{n-4}{q_{n-1}q_{n-2}q_{n-3}}\\
&=& 2n+1+\frac{n-1}{q_{n-1}}+\frac{n-1}{q_{n-1}q_{n-2}}+(n-4)\frac{q_{n-3}-1}
{q_{n-1}q_{n-2}q_{n-3}} \\
&\geq &
2n+1+\frac {n-1}{a_n}+\frac {n-1}{a_n a_{n-1}} +(n-4)
\frac {a_{n-3}-1}{a_n a_{n-1}a_{n-2}}.
\end{eqnarray*}
We would like the right hand side to be at least $(n+2)a_n$. But this is 
equivalent to 
$$(2n+1)a_na_{n-1}a_{n-2}+(n-1)a_{n-1}(a_{n-2}+1)+(n-4)(a_{n-3}+1)-(n+2)
a_n^2 a_{n-1}a_{n-2} \geq 0.$$
Inserting the formulae for $a_n$'s, we get
$$\frac {12(5+\sqrt 5)n^4-2(241+121 \sqrt 5)n^3+2(847+382 \sqrt 5)n^2
-3(341+146 \sqrt 5)n+126+54 \sqrt 5}{(2n-3)(2n-1)(2n+1)(2n+3)^2} \geq 0.$$
The denominator is positive for all integers $n \geq 2$. 
Denote the numerator by $L(n)$ and shift its argument for $6$. The polynomial
$L(n+6)$ has only positive coefficients, so it can not have a positive root.
It then follows that $L(n)$ can not have a root $\gamma \geq 6$. So the
left inequality is valid for all $n\geq 6$, and hence $q_n \geq a_n$.

To prove the other inequality, note that the induction hypothesis implies
$$(n+2)q_n \leq 2n+1+\frac {n-1}{a_{n-1}}+ \frac{n-1}{a_{n-1}}+
\frac{n-1}{a_{n-1}a_{n-2}}+(n-4)\frac {a_{n-2}-1}{a_{n-1}a_{n-2}a_{n-3}}.$$
The condition that the right hand side of this inequality does not exceed 
$(n+2)a_{n+1}$ is equivalent to 
$$(2n+1)a_{n-1}a_{n-2}a_{n-3}+(n-1)a_{n-3}(a_{n-2}+1)+(n-4)(a_{n-2}-1)
-(n+2)a_{n+1}a_{n-1}a_{n-2}a_{n-3} \leq 0.$$
Substituting the formulae for $a_n$'s, we get
$$-3 \frac {(82+42\sqrt 5)n^3-(572+248\sqrt 5)n^2+(1103+474\sqrt 5)n-(529
+247\sqrt 5)}{(2n-3)(2n-1)(2n+1)(2n+5)} \leq 0$$
If we put $n+5$ instead of $n$ in the numerator, we get a polynomial with all
the coefficients positive, and from this we conclude that the numerator does
not change the sign for $n \geq 6$. So, we have proved the inequality $q_n
\leq a_{n+1}$, and thus completed the induction step. This proves the theorem. \qed

{\bf Corollary 3.12}\\ The sequence $q_n=\frac {S^{(1)}(n)}{S^{(1)}(n-1)}$ 
is strictly increasing for all $n \geq 5$, bounded from above 
by $\varphi ^2$ and $\lim_{n \rightarrow \infty} q_n =\varphi ^2=\frac {3+ \sqrt 5}{2}$. \qed

{\bf Remark 3.13}\\ It is proved in \cite{dsv} that the asymptotic behavior
of $S^{(l)} (n)$, $n \geq 0$, is given by $S^{(l)} (n) \sim K_l \alpha _l^n
n^{-3/2}$, where $K_l$ and $\alpha _l$ are constants depending only on the 
rank $l$. Denote $q_n^{(l)}=S^{(l)} (n) / S^{(l)} (n-1)$. Then it follows
$$q^{(l)}_n \sim \alpha _l \left ( 1- \frac {1}{n} \right )^{3/2} \nearrow
\alpha _l.$$
In other words $q_n^{(l)}$ asymptotically behaves as an increasing sequence
tending to $\alpha _l$ as $n \rightarrow \infty$. It can be shown that 
$\alpha _l \in [2,3]$.
The exact values of $\alpha _l$ are known for $l \leq 6$. So, 
$\alpha _0=3$, $\alpha _1=(3+\sqrt 5)/2$, $\alpha _2 =1+\sqrt 2$, etc., and
$\alpha _l \searrow 2$ as $l \rightarrow \infty$.

{\bf Theorem 3.14}\\ The sequences $S^{(l)} (n)$ are log-convex, for $l=2,3$ and $4$. \qed

{\bf Outline of the proof}\\ We present only the case $l = 2$. Dividing the 
short recursion for $S^{(2)} (n)$ by $S^{(2)} (n-1)$ and denoting
the quotient $\frac {S^{(2)} (n)}{S^{(2)} (n-1)}$ by $q_n^{(2)}$, we obtain the 
recursion for the sequence $(q_n^{(2)})$
$$
q_n^{(2)}=\frac{1}{n+2}\left [ 2n+1+\frac{n-1}{q^{(2)}_{n-1}}+\frac{n-4}{q^{(2)}_{n-1}q^{(2)}_{n-2}q^{(2)}_{n-3}}-
\frac{2n-11}{q^{(2)}_{n-1}q^{(2)}_{n-2}q^{(2)}_{n-3}q^{(2)}_{n-4}}-\frac{n-7}{q^{(2)}_{n-1}q^{(2)}_{n-2}q^{(2)}_{n-3}q^{(2)}_{n-4}q^{(2)}_{n-5}}\right ]
$$
with the initial conditions $q^{(2)}_1=q^{(2)}_2=q^{(2)}_3=1$, $q^{(2)}_4=q^{(2)}_5=2$. We want to prove that
the sequence $(q^{(2)}_n)$ is increasing.

From Remark 3.13 we conclude that the sequence $(q^{(2)}_n)$ behaves asymptotically
as $(1+\sqrt 2)(1-\frac{1}{n})^{3/2}$. Denote this quantity by $b_n$, i.e.
$b_n=(1+\sqrt 2)(1-\frac{1}{n})^{3/2}$. Now take the first three terms of the
series expansion of $b_n$ in powers of $\frac{1}{n+2}$. Define
$$a_n=\alpha _2 \left (1-\frac{3}{2} \frac{1}{n+2}+\frac{3}{8}\frac{1}{(n+2)^2}\right )
=\alpha _2 \frac{8n^2+20n+11}{8(n+2)^2},$$
where $\alpha _2 = 1+\sqrt 2$. The sequence $(a_n)$ tends increasingly
toward $\alpha _2$. We shall show now that the sequences $(q^{(2)}_n)$ and $(a_n)$ are
interlaced, i.e. 
$$
a_{n-1} \leq q^{(2)}_n \leq a_n,
$$
for sufficiently large $n$.

Suppose inductively that $a_{n-i-1} \leq q^{(2)}_{n-1} \leq a_{n-i}$ for $i=1,2,3,4,
5$. Then
$$q^{(2)}_n \leq \frac{1}{n+2}\left [ 2n+1+\frac{n-1}{a_{n-2}}+\frac{n-4}{a_{n-2}
a_{n-3}a_{n-4}}-\frac{2n-11}{a_{n-1}a_{n-2}a_{n-3}a_{n-4}}-\frac{n-7}{a_{n-1}
a_{n-2}a_{n-3}a_{n-4}a_{n-5}}\right ].$$
If we prove that the right hand side of this inequality does not exceed $a_n$,
the right inequality, $q^{(2)}_n \leq a_n$,  will follow. But this is equivalent to the condition
$$\frac{P_{10}(n)}{Q_{12}(n)} \geq 0,$$
where $P_{10}(n)$ and $Q_{12}(n)$ are certain polynomials in $n$ of degree
$10$ and $12$, respectively. (Using {\it Mathematica}, the polynomials 
$P_{10}(n)$ and $Q_{12}(n)$ can easily be computed explicitly.) Their leading
coefficients are $262144(1+\sqrt 2)$ and $8^6(1+\sqrt 2)^5$, respectively, and
we can conclude that this quotient is positive for all $n$ big enough. Again,
the biggest real roots of the polynomials $P_{10}(n)$ and $Q_{12}(n)$ can be
easily found using {\it Mathematica}. It turns out that their quotient becomes
(and remains) positive for $n \geq 39$.

On the other hand, from the induction hypothesis and the recursion for $q^{(2)}_n$ it follows that
$$q^{(2)}_n \geq \frac{1}{n+2}\left [ 2n+1+\frac{n-1}{a_{n-1}}+\frac{n-4}{a_{n-1}
a_{n-2}a_{n-3}}-\frac{2n-11}{a_{n-2}a_{n-3}a_{n-4}a_{n-5}}-\frac{n-7}{a_{n-2}a_{n-3}a_{n-4}a_{n-5}a_{n6}}\right ].$$

That the right hand side of this inequality is $\geq a_{n-1}$ is equivalent to
$\frac{P_{13}(n)}{Q_{15}(n)} \geq 0$, where $P_{13}$ and $Q_{15}$ are certain
polynomials in $n$ of degree $13$ and $15$, respectively, with the positive 
leading coefficients. Their quotient is positive for $n$ big enough ($n \geq 6$).

The claim now follows by checking that $(q^{(2)}(n))$ is increasing for $n \leq
44$.

A similar proof works for $l=3$ and $l=4$. We omit the details.
\qed

\section{Calculus method}

\setcounter{equation}{0}
Let us start again as in (3.1) with a linear homogeneous recursion for 
positive numbers $a(n)$, and consider again the corresponding recurrence (3.2)
for the quotients $q(n)=a(n)/a(n-1)$. Suppose again we want to prove that
$a(n)$ is log-convex, i.e. that $(q(n))$ is an increasing sequence (at least 
from some place $n_0$ on).

This time we do the following. Define a continuous function $f : [1,\infty )
\rightarrow \R$ (or $f : [n_0,\infty )\rightarrow \R$) starting with the
appropriate linear function on the segment $[1,2]$ (or $[n_0, n_0+1]$, for
some $n_0 \in \N$) determined by the initial and the next value of $q(n)$, and
then continue to the next segment by the same rule as (3.2). In other words,
by replacing $q \rightarrow f$, $n \rightarrow x$, (3.2) becomes
\be
f(x)=R(x)+\frac{S(x)}{f(x-1)},
\ee
defined so for $x \geq 2$ (or $x \geq n_0 +1$). In a sense, $f$ is a dynamical
system patching the discrete values $f(n)=q(n)$ if $R$ and $S$ are ``good''
enough functions. For example, if $R$ and $S$ are rational functions (and in
combinatorics this is mostly the case) without poles on the positive axis, then
$f$ is a piecewise rational function, i.e. rational on every open interval
$(n,n+1)$ for any integer $n \geq 1$ (or $n \geq n_0$). If we can prove that
$f$ is smooth on such open intervals (usually by proving that $f$ is
bounded from above and below by some well-behaved functions), we can consider 
the derivative $f'(x)$, for any $x \in (n,n+1)$. The idea is to show inductively
on $n$ that $f'(x) \geq 0$ (or $f'(x) \leq 0$ if we want to prove log-concavity
of $a(n)$). This will imply that $f$ is an increasing function (or a decreasing
function in the log-concave case) on any open interval $(n,n+1)$, and then, by
continuity of $f$ it will follow that $f$ is increasing (or decreasing) on 
its whole domain, and hence $q(n)=f(n)$, $n \in \N$, will increase (or decrease), too.

Now, in general, if we want to prove that a sequence $q(n)$ defined by (3.2) is
increasing, we form the corresponding functional equation (4.1) with the
appropriate start, i.e. we define $f(x)$ on some starting segment $[n_0,n_0+1]$
to be an increasing function and then prove inductively on $n$ that $f'(x)\geq
0$ for $x\in (n,n+1)$. It is always necessary to have some a priori bounds
for $f$. So, suppose we know $0 < m(x) \leq f(x) \leq M(x)$, for all $x \geq
n_0$. Let us find some sufficient conditions which ensure that $f'(x) \geq 0$.
Of course, we assume that $R$ and $S$ are smooth on all open intervals
$(n,n+1)$, $n \geq n_0$.

Fix an $x \in (n,n+1)$ and write $f=f(x)$, $f_1=f(x-1)$, $R=R(x)$, $S=S(x)$.
Then (4.1) can be written as $$f=R+\frac{S}{f_1},$$ or equivalently
\be
ff_1=Rf_1+S.
\ee
Taking the derivative $d/dx$ of both sides of (4.2), we get
$$f'f_1+ff_1'=R'f_1+Rf_1'+S',$$ implying $$f'f_1=R'f_1+S'+(R-f)f_1'.$$
From (4.2), we obtain
\be
f'f_1=R'f_1+S'-\frac{S}{f_1}f_1'.
\ee
Assume that 
\be0 < m(x) \leq f(x) \leq M(x),
\ee
for all $x \geq n_0$, and suppose inductively that $f$ is increasing on
$[n_0, n]$. This means that $f_1' \geq 0$. Then, with our notation convention
$m_1=m(x-1)$ etc., we have:

{\bf Theorem 4.1}\\ If $R' \geq 0$, $R'm_1+S' \geq 0$ and $S \leq 0$, then 
$f' \geq 0$.

{\bf Proof}\\ Obvious from (4.3), (4.4) and the above discussion. \qed

{\bf Theorem 4.2}\\ Suppose $R' \geq 0$, $S' \geq 0$, $S \geq 0$ and $m_1m_2
(R'm_1+S')\geq S(R_1'M_2+S_1')$. Then $f' \geq 0$.

{\bf Proof}\\ Divide (4.3) by $f_1$ and substitute in by the same rule $f_1'$.
WE obtain
$$f'=\frac{R'f_1+S'}{f_1}-\frac{S}{f_1^2}f_1'=\frac{R'f_1+S'}{f_1}-\frac{S}{f_1^2}
\left [ \frac{R_1'f_2+S_1'}{f_2}-\frac{S_1}{f_1^2}f_2' \right ],$$
and this implies
\be
f'=\frac{R'f_1+S'}{f_1}-\frac{S}{f_1^2}\frac{R_1'f_2+S_1'}{f_2}+\frac{SS_1}
{(f_1f_2)^2}f_2'.
\ee
By the inductive hypothesis, $f_2'=f'(x-2)\geq 0$, and since $SS_1 \geq 0$, it
follows that the last term in (4.5) is non-negative. On the other hand,
\be
\frac{R'f_1+S'}{f_1}-\frac{S}{f_1^2}\frac{R_1'f_2+S_1'}{f_2} \geq 0
\Longleftrightarrow f_1f_2(R'f_1+S') \geq S(R_1'f_2+S_1').
\ee
Since $R'$, $S'$, $S \geq 0$ and $m_1m_2(R'm_1+S')\geq S(R_1'M_2+S_1')$, by our
assumptions, it follows that (4.6) holds, and we are done. \qed

Now let us show how to apply the above method to some combinatorially
relevant numbers.

Recall that the $n$-th {\bf big Schr\"oder number} $r_n$ is the number of
lattice paths from $(0,0)$ to $(n,n)$ with steps $(1,0)$, $(0,1)$ and
$(1,1)$ that never rise above the line $y=x$. Equivalently, $r_n$ is the 
number of lattice paths from $(0,0)$ to $(2n,0)$ with steps $(1,1)$,
$(1,-1)$ and $(2,0)$ that never fall below the $x$-axis. For $n=0$ we put $r_0
:=1$. So, $r_n$ is the number of Motzkin paths whose every plateau, every
valley, every terrace and every plain is of even length. Hence, via Proposition
3.6, big Schr\"oder numbers count secondary structures whose unpaired bases 
form contiguous blocks of even length. The $n$-th {\bf little Schr\"oder 
number} is $s_n=\frac{1}{2}r_n$, $s_0=1$.

{\bf Theorem 4.3}\\ The Schr\"oder numbers are log-convex. The sequences
$r_n/r_{n-1}$ and $s_n/s_{n-1}$ increasingly tend to $3+2 \sqrt 2$ as $n
\rightarrow \infty$.

{\bf Proof}\\ Similarly to the proof of Proposition 3.2, it can be shown (and
it is well known) that $(r_n)$ satisfies the following (convolutional) 
recurrence:
$$r_{n+1}=r_n+\sum _{j=0}^n r_jr_{n-j}.$$
This implies that the generating function for $(r_n)$ is given by
$$r(x)=\frac{1-x-\sqrt{1-6x+x^2}}{2x}.$$
This, in turn, implies (as in Proposition 3.2) that $r_n$ satisfy the following
short recursion
\be
r_n=\frac{1}{n+1}[3(2n-1)r_{n-1}-(n-2)r_{n-2}],
\ee
together with the initial conditions $r_0=1$, $r_1=2$. (There are also
combinatorial proofs of (4.7), see \cite{foata} or \cite{shapsul}.) 
Divide (4.7) by $r_{n-1}$
and denote $r_n/r_{n-1}$ by $q(n)$. We get
\be
q(n)=\frac{3(2n-1)}{n+1}-\frac{n-2}{n+1}\frac{1}{q(n-1)}, \quad n \geq 2,
\ee
with $q(1)=2$. Then $q(2)=3$ and $q(3)=11/3$. Define the function $f : [2,
\infty) \rightarrow \R$ by
$$f(x)= \left\{\matrix{
\hfill \frac {1}{3} (2x+5) &,& x \in [2,3] \cr
\frac{1}{x+1} [6x-3-\frac {x-2}{f(x-1)}] &,& \hfill x\geq 3. \cr
} \right. $$
It is easy to show by induction on $n$ that $f$ is bounded on $[2,n]$, for $n
\geq 2$. More precisely, $3 \leq f(x) \leq 6$, for all $x \geq 2$. Also, $f$ is
continuous everywhere and differentiable on open intervals $(n,n+1)$, $n \geq 2$. In
the notations of Theorem 4.1, we have
$$R(x)=\frac{3(2x-1)}{x+1}, \quad S(x)=-\frac{x-2}{x+1}, \quad m(x)=3, \quad M(x)=6.$$
Then
$$R'(x)=\frac{9}{(x+1)^2} \geq 0, \quad S'(x)=-\frac{3}{(x+1)^2}.$$
Clearly, $S(x)\leq 0$ for $x \geq 2$, and let us check that $R'm_1+S' \geq 0$.
But this is obvious, since
$$\frac{9}{(x+1)^2} \cdot 3 -\frac{3}{(x+1)^2} =\frac{24}{(x+1)^2} \geq 0.$$
Since $f'(x)=2/3 \geq 0$ for $x \in (2,3)$, it follows from Theorem 4.1 that
$f'(x) \geq 0$ for all $x \in (n,n+1)$, $n \geq 2$. Hence, by continuity, $f$
is increasing, so $(q(n))$ is increasing and therefore $(r_n)$ is log-convex.
By passing to limit in (4.9) we get the last claim. \qed

Now that we are more familiar with this method, let us show once again (but
this time almost automatically) the log-convexity of Motzkin numbers.

Recalling the short recursion (3.9), we define 
the function $f: [2, \infty ) \rightarrow \R$ by
$$f(x)= \left\{\matrix{
\hfill 2 &,& x \in [2,3] \cr
\frac{1}{x+2} [2x+1+\frac {3(x-1)}{f(x-1)}] &,& \hfill x\geq 3 .\cr
} \right. $$
Here we have
$$R(x)=\frac{2x+1}{x+2}, \quad S(x)=\frac{3(x-1)}{x+2}, $$
so
$$R'(x)=\frac{3}{(x+2)^2} \geq 0, \quad S'(x)=\frac{9}{(x+2)^2} \geq 0.$$
Further,it is easy to check that $2 \leq f(x) \leq 7/2$ for all $x \geq 2$.
So, we may take $m(x)=2$, $M(x)=7/2$. Let us check the inequality $m_1m_2(R'm_1+S')
\geq S(R_1'M_2+S_1')$ from Theorem 4.2. But this is equivalent to
$$4\left [ 2\cdot \frac{3}{(x+2)^2}+\frac{9}{(x+2)^2} \right ] \geq
\frac{3x-3}{x+2}\left [ \frac{7}{2} \frac{3}{(x+2)^2}+\frac{9}{(x+2)^2} \right ],$$
and this is equivalent to $$1.5 x^2+61.5 x+117 \geq 0.$$
This last inequality is certainly true for all $x \geq 1$. Hence, by Theorem 
4.2, $f' \geq 0$, and this implies the log-convexity of Motzkin numbers. \qed

A {\bf directed animal} of size $n$ is a subset $S \subseteq (\N \cup \{ 0 \} )
^2$ of cardinality $n$ with the following property: if $p \in S$, then there is
a lattice path from $(0,0)$ to $p$ with steps $(1,0)$, $(0,1)$, all of whose
vertices lie in $S$. Let $a_n$ be the number of directed animals of size $n$.
It is well known that the generating function for $(a_n)_{n \geq 1}$ is given
by
$$A(x)=\sum _{n \geq 1} a_nx^n=\frac{1}{2}\left( \sqrt{\frac{1+x}{1-3x}}-1
\right ) = x+2x^2+5x^3+13x^4+35x^5+\ldots $$
These same numbers appear also as the row sums of Motzkin triangle (see \cite
{donaghey} and \cite{barcucci}).

{\bf Theorem 4.4}\\ The sequence $(a_n)_{n \geq 1}$ of numbers of directed
animals is log-convex. Hence the sequence $a_n/a_{n-1}$ increasingly tends to
$3$.

{\bf Proof}\\ Taking the derivative $A'(x)$ of the generating function $A(x)$
we get $$(1+x)(1-3x)A'(x)=2A(x)+1,$$ and equating coefficients of $x^{n-1}$
above yields the recurrence
$$na_n=2na_{n-1}+3(n-2)a_{n-2}, \quad n \geq 3,$$
with $a_1=1$, $a_2=2$.

The corresponding function for the successive quotients is given by
$$f(x)=2+\frac{3(x-2)}{x}\frac{1}{f(x-1)}, \quad x \geq 3,$$
and $f(x)=\frac{x+2}{2}$ on $[2,3]$. It is easy to check that $2 \leq f(x) \leq
\frac{7}{2}$.

Here $R(x)=2$, $S(x)=\frac{3(x-2)}{x} \geq 0$ and $R'(x)=0$, $S'(x)=\frac{6}{x^2} \geq 0$, and $m=2$, $M=7/2$, so the claim follows from Theorem 4.2. \qed

We aplly now our approach to Franel's sequences. Recall that the $n$-th 
{\bf Franel number} of order $r$ is defined by
$$S_n^{(r)}=\sum _{k=0}^n \binom{n}{k}^r.$$
The numbers $S_n^{(r)}$ satisfy a linear homogenous recurrence of order $\lfloor
\frac{r+1}{2} \rfloor$ with polynomial coefficients. For $r=3$, the recurrence
reads as follows (\cite{stanleyII}):
\be
n^2S_n^{(3)}=(7n^2-7n+2)S_{n-1}^{(3)}+8(n-1)^2S_{n-2}^{(3)},
\ee
with the initial conditions $S_0^{(3)}=1$, $S_1^{(3)}=2$.

{\bf Theorem 4.5}\\ The sequence $\left ( S_n^{(3)} \right ) _{n \geq 0}$ is
log-convex. The same is true for the sequence $\left ( S_n^{(4)} \right ) _{n \geq 0}$.

{\bf Proof}\\ Again, taking the quotients $q(n)=S_n^{(3)}/S_{n-1}^{(3)}$ in
(4.9), we get
\be
q(n)=\frac{7n^2-7n+2}{n^2}+\frac{8(n-1)^2}{n^2}\frac{1}{q(n-1)}, \quad n \geq 2,
\ee
with $q(1)=2$. Then $q(2)=5$. Define the function 
$f: [1, \infty ) \rightarrow \R$ by
$$f(x)= \left\{\matrix{
\hfill 3x-1 &,& x \in [1,2] \cr
\frac{7x^2-7x+2}{x^2} +\frac {8(x-1)^2}{x^2}\frac{1}{f(x-1)}] &,& \hfill x\geq 2. \cr
} \right. $$
The function $f$ is bounded. More precisely, $5 \leq f(x) \leq 9$ for all $x
\geq 2$. Indeed, $f$ is clearly bounded on $[1,2]$. On $[2,3]$ $f$ is given by
$$f(x)=\frac{21x^2-41x+18}{3x^2-4x},$$
and one can check easily that $5 \leq f(x) \leq 9$ on this interval. Now, by
induction on $n \geq 3$, it is not hard to verify that if $5 \leq f(x) \leq 9$
for $x \in [2,n]$, then also $5 \leq f(x) \leq 9$ for $x \in [n,n+1]$. Also, 
$f$ is a continuous function. So, here we have, in the notations of Theorem 4.1:
$$R(x)=\frac{7x^2-7x+2}{x^2}, \quad S(x)=\frac{8(x-1)^2}{x^2}$$
and we may take $m(x)=5$, $M(x)=9$ for $x\geq 2$. Since $R'(x)=\frac{7x-2}{x^3}
\geq 0$, $S'(x)=\frac{16(x-1)}{x^3} \geq 0$ and $S(x) \geq 0$ for $x \geq 2$,
it remains only to check the last condition in Theorem 4.2, i.e.
$$25\left [ 5\cdot \frac{7x-2}{x^3}+\frac{16(x-1)}{x^3} \right ] \geq
\frac{8(x-1)^2}{x^2}\left [ \frac{7x-9}{(x-1)^3} \cdot 9 +\frac{16(x-2)}{(x-1)^3} \right ].$$
This is equivalent to $$643x^2-1021x+650 \geq 0,$$
and this last inequality is certainly true for $x \geq 1$.

So, as before, we conclude that $f$ increases, hence $\left ( q(n) \right ) 
_{n \geq 0}$ is an increasing sequence and therefore $\left ( S_n^{(3)} \right ) _{n \geq 0}$ is log-convex.

Note that $(q(n))$ converges, being increasing and bounded, and so, passing to
limit in (4.10) we obtain that
$$\lim _{n \rightarrow \infty} \frac{S_n^{(3)}}{S_{n-1}^{(3)}} = 8.$$

A similar proof applies to the case of $\left ( S_n^{(4)} \right ) _{n \geq 0}$,
starting with the recurrence
$$n^3S_n^{(4)}=2[6n^3-9n^2+5n-10]S_{n-1}^{(4)}+(4n-3)(4n-4)(4n-5)S_{n-2}^{(4)},$$
$S_0^{(4)}=1$, $S_1^{(4)}=2$. It turns out that
$$\lim _{n \rightarrow \infty} \frac{S_n^{(4)}}{S_{n-1}^{(4)}} = 16.$$
We leave the details to the reader. \qed

As a bit more simple case, let us show the log-behavior of sequences defined
by two-term recurrences with constant coefficients.

\newpage
{\bf Proposition 4.6}\\ Let $\left ( a_n \right ) _{n \geq 0}$ be a positive
sequence defined by $$ a_n=C_1 a_{n-1} -C_2 a_{n-2},$$ where $C_1$, $C_2 > 0$ 
are constants, and with some initial conditions $a_0 > 0$, $a_1 > 0$. Then the log-
behavior of $(a_n)$ is completely determined by the log-behavior of its first
three terms. In other words, the sequence $(a_n)$ is log-convex if $a_1^2 \leq
a_0a_2$, and log-concave if $a_1^2 \geq a_0a_2$.

{\bf Proof}\\ Let $q(n)=a_n /a_{n-1}$. The recursion for $q(n)$ is given by
$$q(n)=C_1-\frac{C_2}{q(n-1)},$$ and the first two values are $q(1)=a_1/a_0$,
$q(2)=a_2/a_1$. Set $\Delta = a_0a_2-a_1^2$ and define the function
$f: [1, \infty ) \rightarrow \R$ by
$$f(x)= \left\{\matrix{
\hfill \frac{\Delta }{a_0a_1}(x-1)+\frac{a_1}{a_0} &,& x \in [1,2] \cr
C_1-\frac {C_2}{f(x-1)} &,& \hfill x\geq 2 .\cr
} \right. $$
$f$ is a linear fractional mapping on any interval $(n,n+1)$, $n \in \N$, and
$f$ is continuous.

Suppose first that $\Delta \geq 0$. Then $f$ is obviously increasing on $[1,2]$
and bounded from below by $f(1)=q(1)=a_1/a_0$. Namely, suppose that $f(x) \geq
f(1)$ on $[1,n]$. For $x \in [n,n+1]$ we have
$$f(x) =C_1-\frac {C_2}{f(x-1)} \geq C_1-\frac {C_2}{f(1)} = f(2) \geq f(1),$$
and so $f(x)$ is bounded from below on the whole interval $[1,n+1]$. Hence $f$
is bounded from below by $f(1)$. This also shows that $f$ is differentiable on 
any open interval $(n,n+1)$, $n \in \N$. Here we have $R(x)=C_1$, $S(x)=-C_2$
(so $R'(x)=S'(x)=0$ for $x \not\in \N$) and we may take $m(x)=a_1/a_0$. It
follows by Theorem 4.1 that $f'(x) \geq 0$ for all $x \not\in \N$. By 
continuity, $f$ is an increasing function, and hence $f(n)=q(n)$ is an 
increasing sequence.

If $\Delta \leq 0$, then $f$ decreases on $[1,2]$ and similarly as above $f$ is,
in this case, bounded from above by $f(1)$. Since $$f'(x)=C_2\frac{f'(x-1)}
{f(x-1)},$$ we see that the sign of the derivative is propagated to the right
and is determined by its sign on $(1,2)$. Hence, in this case, $f'(x) < 0$. \qed

{\bf Corollary 4.7}\\ The sequence $\left ( F_{2n+1} \right ) _{n \geq 0}$ of
odd-indexed Fibonacci numbers is log-convex, while the sequence of even-indexed
Fibonacci numbers $\left ( F_{2n} \right ) _{n \geq 1}$ is log-concave.

{\bf Proof}\\ Let $a_n=F_{2n+1}$. The numbers $a_n$ satisfy the recursion
$a_n = 3a_{n-1}-a_{n-2}$, and $a_1^2=4<5=a_0a_2$. The numbers $b_n=F_{2n}$
satisfy the same recursion $b_n = 3b_{n-1}-b_{n-2}$, but this time with 
initial conditions $b_1=1$, $b_2=3$. However, $b_2^2=9>8=b_1b_3$.
So, the claim follows from Proposition 4.6. (Of course, this Corollary is an
immediate consequence of the Cassini identity $F_n^2-F_{n-1}F_{n+1}=(-1)^n$.)
\qed

As our next application of the method, we consider some classical orthogonal
polynomials. Let $\nu > -1/2$ be a parameter. The {\bf Gegenbauer} (or
{\bf ultraspherical}) {\bf polynomials} $C_n^{(\nu )}(t)$ are defined by 
(\cite{szego})
\be
nC_n^{(\nu )}(t)=2t(\nu +n-1)C_{n-1}^{(\nu )}(t)-(2\nu +n-2)C_{n-2}^{(\nu )}(t), \quad n\geq 2,
\ee
together with the initial conditions $C_0^{(\nu )}(t)=1$, $C_1^{(\nu )}(t)=
2\nu t$. 

For $\nu = 1/2$, $C_n^{(1/2)}(t)=P_n(t)$ is the {\bf Legendre
polynomial} and for $\nu = 1$, $C_n^{(1)}(t)=U_n(t)$ is the {\bf Chebyshev
polynomial of the second kind}.

\newpage
{\bf Theorem 4.8}\\ The sequence $\left ( C^{(\nu )} _n(t)\right )_{n \geq 0}$ is
log-concave for $t \geq 1$, $\nu \geq 1$, and log-convex for $0< \nu < 1$,
$t \geq \max \left \{ \frac{1}{\sqrt{2\nu }},\frac{1}{\sqrt{2(1-\nu )}} \right \}$.

{\bf Proof}\\ Let $q_t(n)= \frac{C^{(\nu )} _n(t)}{C
^{(\nu )} _{n-1}(t)}$. 
Dividing the recursion (4.11) by $n C^{(\nu )} _{n-1} (t)$ we obtain the
recursion
\be
q_t(n)=2t\frac{\nu+n-1}{n}-\frac{2\nu +n-2}{n}\frac{1}{q_t(n-1)}, n \geq 1
\ee
with initial condition $q_t(1)=2\nu t$.
Then $q_t(2)=t(\nu +1) - \frac{1}{2t}$.
Now define the function $f_t:[1,\infty ) \rightarrow \R$ by formula
\be
f_t(x)= \cases {
\frac{2t^2(1-\nu )-1}{2t}x+\frac{2t^2(3\nu -1)+1}{2t} &, if $x\in [1,2]$,\cr
\frac{2t(\nu+x-1)}{x} - \frac {2\nu +x-2}{x}\frac {1}{f_t(x-1)} &, if $x\geq 2.$\cr
}
\ee
The function $f_t(x)$ is piecewise rational, i.e. it is rational on all segments
$[n,n+1]$, $n \in \N$. It is easy to show, by induction on $n$, that $f_t(x)$ is
continuous at the points $x=n$, for all $n \in \N$, and that $f_t(n)=q_t(n)$.

Let us first consider the case $\nu \geq 1$ and $t \geq 1$. We claim that $f_t(x)$ is smooth
on all intervals $(n,n+1)$, $n\in \N$. To prove this, it is enough to show that
$f_t(x)$ is bounded from below by some positive quantity. Let us check that
$f_t(x) \geq \frac{1}{t}$ for all $x \geq 1$.

Since $\nu \geq 1$, $t \geq 1$, the function $f_t(x)$ is non-increasing on the interval
$[1,2]$, so it is enough to show that $f_t(2) \geq \frac{1}{t}$. But this is
equivalent to $t(\nu +1) > \frac{3}{2t}$, or $2t^2(\nu +1) \geq 3$. Since $\nu \geq 1$,
we have $4t^2 \geq 3$, and this is true for all $t \geq 1$.

Suppose now inductively that $f_t(x) \geq \frac{1}{t}$ for all $x \in [1,n]$,
and take an $x\in [n,n+1]$, $n \geq 2$. From (4.13) we have
$$f_t(x) \geq \frac{2t(\nu+x-1)}{x} - \frac {2\nu +x-2}{x}\frac {1}{\frac{1}{t}}=t,$$
and $t\geq \frac{1}{t}$ for all $t \geq 1$.
So, the function $f_t(x)$ is nowhere zero on $[1,\infty )$, hence it has no 
poles on $[1,\infty )$ and hence is smooth on all intervals $(n,n+1)$, $n\in \N$.

The function $f_t(x)$ is obviously decreasing on $[1,2]$. Suppose now that it
is decreasing on $[1,n]$. Let $x \in (n,n+1)$. Then
\be
f_t'(x)=\frac{2(1-\nu )}{x^2f_t(x-1)}\left [ t f_t(x-1)-1\right ] +
\left (1+2\frac{\nu -1}{x}\right )\frac {f_t'(x-1)}{f_t^2(x-1)}.
\ee
The term in square brackets is positive, since $f_t(x) \geq \frac{1}{t}$, so the
whole first term in (4.14) is negative. The second term is negative by the induction
hypothesis, hence $f_t'(x) \leq 0$. This completes the step of induction, and we
can conclude that the function $f_t(x)$ is decreasing on $[1,\infty )$. Hence,
the sequence $\left ( C^{(\nu )} _n(t) \right )_{n \geq 0}$ is log-concave for all
$\nu \geq 1$, $t \geq 1$.

Let us now consider the case $0 < \nu < 1$,
$t \geq \max \left \{ \frac{1}{\sqrt{2\nu }},\frac{1}{\sqrt{2(1-\nu )}} \right \}$.
Obviously, for such values of $t$, the function $f_t(x)$ is increasing on
$[1,2]$, and $f_t(1) \geq \frac{1}{t}$, hence $f_t(x) \geq \frac{1}{t}$ on
$[1,2]$. It follows easily by induction on $n$ that $f_t(x) \geq \frac{1}{t}$
on $[1,\infty )$, hence $f_t(x)$ is bounded from below by a positive quantity,
hence $f_t(x)$ is smooth on all open intervals $(n,n+1)$, $n\in \N$.

Suppose now that $f_t(x)$ is increasing on $[1,n]$ and let $x \in (n,n+1)$, $
n \geq 2$.
Take a look at (4.14) again.
The first term is positive since $f_t(x) \geq \frac{1}{t}$, and the second term is
positive by the induction hypothesis for all $x\geq 2$. Hence $f_t'(x)\geq 0$
and $f_t(x)$ is increasing on $[1,n+1]$. The increasing behavior of $f_t(x)$
implies the log-convexity of the sequence $\left ( C^{(\nu )} _n(t) \right )_{n \geq 0}$. \qed

In Fig. 5 we can see the areas of log-convexity (denoted by $\Lambda $) and
log-concavity ($\Xi $) of sequences $\left ( C^{(\nu )} _n(t) \right )_{n \geq 0}$.
Our analysis is not sufficient to determine the logarithmic behavior of 
sequences $\left ( C^{(\nu )} _n(t) \right )_{n \geq 0}$ whose parameters $(\nu, t)$
fall into the areas denoted by I and II.

\vskip 5mm
\begin{figure}[h] \centerline {
\epsfig{file=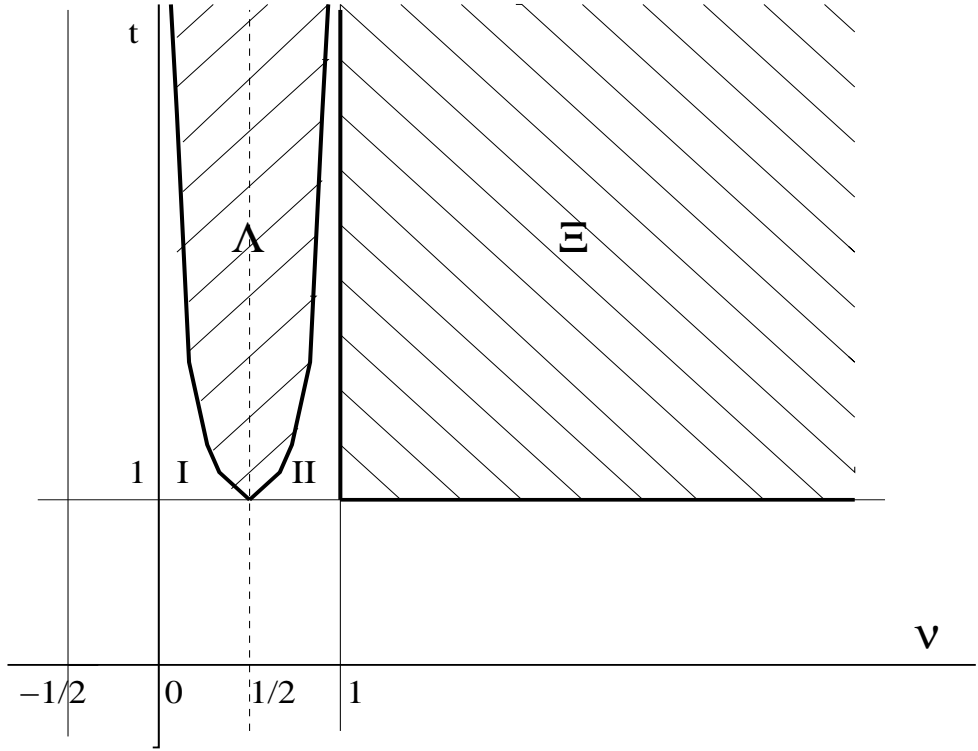,height=10cm,width=13cm,silent=}}
\caption{Areas of log-convexity ($\Lambda $) and log-concavity ($\Xi $) for
sequences of values of Gegenbauer polynomials}
\label{gegen}
\end{figure}

{\bf Remark 4.9}\\ So far we have restricted our definitions of log-convexity 
and log-concavity to positive (or perhaps non-negative) sequences, since most
combinatorial sequences have only non-negative terms. The definitions,
however, make sense even if we include negative numbers.

{\bf Corollary 4.10}\\ The sequence $\left ( C^{(\nu )} _n(t)\right )_{n \geq 0}$ is
log-concave for $|t| \geq 1$, $\nu \geq 1$, and log-convex for $0< \nu < 1$,
$|t| \geq \max \left \{ \frac{1}{\sqrt{2\nu }},\frac{1}{\sqrt{2(1-\nu )}} \right \}$.
\newpage
{\bf Proof} \\ It follows from Theorem 4.7 and the relation
$$C^{(\nu )} _n(-t) = (-1)^nC^{(\nu )} _n(t),$$
valid for all $\nu \neq 0$ and all $t$. \qed

{\bf Corollary 4.11}\\ The sequence $\left ( U_n(t)\right )_{n \geq 0}$ of the
values of Chebyshev polynomials of the second kind is log-concave for $|t| 
\geq 1$. In particular, $F_{2n}=U_n(3/2)$ is log-concave. 

{\bf Proof} \\ Follows from Corollary 4.10, since $U_n(t)=C^{(1)} _n(t)$. \qed

{\bf Corollary 4.12}\\ The sequence $\left ( P_n(t)\right )_{n \geq 0}$ of the
values of Legendre polynomials is log-convex for $|t| \geq 1$ and 
log-concave for $|t| \leq 1$.

{\bf Proof} \\ $P_n(t)=C^{(1/2)} _n(t)$ and log-convexity follows from
Corollary 4.10. The log-concavity follows from the well known fact that
$|P_n(t)| \leq 1$ for $|t| \leq 1$ and from the inequality
(\cite{abramowitz}, \cite{szasz})
$$P_n^2(t)-P_{n-1}(t)P_{n+1}(t) \geq \frac{1-P_n^2(t)}{(2n-1)(n+1)},$$
valid for $|t| \leq 1$. \qed

Our results on Legendre polynomials complement and generalize in some respects those of
Tur\'{a}n and Szeg\"o (\cite{szego1}), who considered determinants of the form $\left |
{{P_n(x)} \atop {P_{n-1}(x)}} \quad {{P_{n+1}(x)} \atop {P_n(x)}} \right |$.

Now, a combinatorial consequence is in order. 
The $n$-th {\bf central Delannoy number} $D(n)$ counts the lattice
paths in $(x,y)$ coordinate plane from $(0,0)$ to $(n,n)$ with steps $(1,0)$, 
$(0,1)$ and $(1,1)$. (Such paths are also known as {\bf ``king's paths''}.)
The generating function of $\left ( D(n) \right ) _{n \geq 0}$ is given by
$D(x)=\frac {1}{\sqrt {1-6x+x^2}}$ (\cite{stanleyII}). This implies easily that
$D(n)$'s satisfy the recurrence (4.11) for $\nu =1/2$ and $t=3$ and $D(0)=1$,
$D(1)=3$. Hence $D(n)=P_n(3)$, the value of the $n$-th Legendre polynomial
$P_n$ at $t=3$.

{\bf Corollary 4.13.}\\ The sequence $D(n)$ of Delannoy numbers is log-convex
and 
$\lim_{n \rightarrow \infty} \frac {D(n)}{D(n-1)}=3+2 \sqrt 2$.

{\bf Proof}\\ The log-convexity follows from Corollary 4.12. Hence $q_3(n)=
D(n)/D(n-1)$ is increasing. From (4.12) for $\nu =1/2$ and $t=3$ it follows
easily by induction that $q_3(n)\leq 6$. Passing to limit, we obtain the
second claim. \qed

{\bf Remark 4.14}\\ No wonder that the above limit coincides with the limit
for Schr\"oder numbers from Theorem 4.3. This can be explained by showing that
$$r_n=\frac{P_n'(3)}{P_n'(1)}.$$
(For a proof, see e.g. \cite{doslicphd}.) 

The derivatives of ultraspherical polynomials behave log-concavely. More
precisely, the following is valid.

{\bf Theorem 4.15}\\ The sequence $\left ( C^{(\nu )'} _n(t) \right )_{n\geq 0}$
is log-concave, for all $\nu >0$, $t \geq 1$.

{\bf Proof}\\
Deriving the recurrence (4.11) with respect to $t$ and using the relation
$$2(\nu +n-1) C^{(\nu )}_{n-1}(t)=C^{(\nu )'}_n(t)-C^{(\nu )'}_{n-2}(t),$$
we get the two-term recursion for derivatives:
$$C^{(\nu )'} _n(t)=2t\frac{\nu +n-1}{n-1}C^{(\nu )'} _{n-1}(t)-\frac{2\nu +n-1}{n-1}C^{(\nu )'} _{n-2}(t),$$
which starts with $C^{(\nu )'}_0(t)=0$, $C^{(\nu )'}_1(t)=2 \nu $, 
$C^{(\nu )'}_2(t)=4\nu(\nu +1)t$ and $C^{(\nu )'}_3(t)=4\nu (\nu +1)(\nu +2)t^2-2\nu (\nu +1).$ Hence, the successive quotients $q_t(n)=C^{(\nu )'} _n(t)/
C^{(\nu )'} _{n-1}(t)$ satisfy the recurrence
$$q_t(n)=2t\frac{\nu +n-1}{n-1}-\frac{2\nu +n-1}{n-1}\frac{1}{q_t(n-1)}, n\geq 3,$$
starting with $q_t(2)=2(\nu +1)t$, $q_t(3)=(\nu +2)t-\frac{1}{2t}$.
The appropriate function $f_t:[2,\infty)\rightarrow \R$ 
is defined by
$$ 
f_t(x)= \cases {
-\frac{1}{2t}(2\nu t^2+1)x+(4\nu +1)t+\frac{1}{t} &, if $x\in [2,3]$,\cr
\frac{2t(\nu +x-1)}{x-1} - \frac {2\nu +x-1}{x-1}\frac {1}{f_t(x-1)} &, if $x\geq 3.$\cr
}
$$ 
The function $f_t(x)$ is obviously decreasing on $[2,3]$ for all $\nu > 0$, $t
\geq 1$. Also, $f_t(x)$ is continuous in all points $x=n$, for $n\in \N$, and
$f_t(n)=q_t(n)$. As before, it is not hard to check that $f_t(x) \geq \frac{1}{t}$, for all $x\geq 2$.
Hence, the derivative $f_t'(x)$ exists for any $x \in (n,n+1)$, $n \geq 2$.
$$f_t'(x)=\frac{2\nu }{(x-1)^2f_t(x-1)}\left [ 1-t f_t(x-1)\right ]
+\left ( 1+\frac{2\nu }{x-1}\right )\frac{f_t'(x-1)}{f_t^2(x-1)}.$$
This derivative is obviously non-positive on $(2,3)$. Suppose 
inductively that $f_t'(x-1) \leq 0$.
The first term above is negative since $f_t(x) \geq \frac{1}{t}$, and the second term
is negative by the induction hypothesis. Hence, $f_t(x)$ is decreasing on 
$(n,n+1)$, and then, by continuity, on $[1,n+1]$. This completes the inductive
step and hence $q_t(n)=f_t(n)$ is also decreasing. \qed

We conclude our review of orthogonal polynomials with
Laguerre polynomials.
The (ordinary) {\bf Laguerre polynomials} satisfy the recursion
$$n L_n(t)=(2n-1-t)L_{n-1}(t)-(n-1)L_{n-2}(t)$$
with the initial conditions $L_0(t)=1$, $L_1(t)=1-t$.

{\bf Theorem 4.16}\\ The sequence $\left (L_n(t)\right )_{n\geq 0}$ is 
log-concave for all $t \leq 0$. The ratio $L_n(t)/L_{n-1}(t)$ is greater than
one for all $t < 0$ and it tends to $1$ decreasingly as $n \rightarrow \infty$.

{\bf Proof}\\
Dividing the above recursion by $nL_{n-1}(t)$ and denoting $\frac{L_n(t)}{L_{n-1}(t)}$
by $q_t(n)$, we get a recursion for $q_t(n)$
\be 
q_t(n)=\frac{2n-1-t}{n}-\frac{n-1}{n}\frac{1}{q_t(n-1)}
\ee
with $q_t(1)=1-t$. By computing $q_t(2)=\frac{t^2-4t+2}{2(1-t)}$, we see that
$q_t(1) \geq q_t(2)$, for all $t \leq 0$.
Define $f_t:[1,\infty) \rightarrow \R$ by formula
$$f_t(x)= \cases {
-\frac{t^2}{2(t-1)}x+\frac{3t^2-4t+2}{2(t-1)} &, if $x\in [1,2]$,\cr
\frac{2x-1-t}{x} - \frac {x-1}{x}\frac {1}{f_t(x-1)} &, if $x\geq 2$.\cr
}$$
This function is, again, piecewise rational, continuous for $x=n \in \N$
and $f_t(n)=q_t(n)$ for all $n \in \N$. By induction on $n \geq 2$, one can
easily check that, for $x\in [1,n]$, we have
$1 \leq f_t(x) \leq 1-t$ for all $t \leq 0$. Hence,
\be
1 \leq f_t(x) \leq 1-t, \quad x\geq 1, t\leq 0.
\ee
By computing the derivative, we see that $f_t'(x) \leq 0$ for $x\in (1,2)$, 
and for $x \in (n,n+1)$, $n \geq 2$ we find
$$f_t'(x)=\frac{1+t}{x^2}-\frac{1}{x^2}\frac{1}{f_t(x-1)}+\left ( 1-\frac{1}{x}
\right )\frac{f_t'(x-1)}{f_t^2(x-1)}.$$
Suppose $f_t'(x-1) \leq 0$. Then for $t\leq -1$, all three 
terms on the right hand side are negative, hence
$f_t'(x) \leq 0$. For $-1 < t \leq 0$, the second and third terms are still
negative, but the first term is positive. The claim will follow if we prove
$$\frac{1+t}{x^2}-\frac{1}{x^2}\frac{1}{f_t(x-1)} \leq 0.$$
But this is equivalent to $(1+t)f_t(x-1)-1 \leq 0$, and this reduces to 
$f_t(x-1) \leq \frac{1}{1+t}$. But $f_t(x) \leq 1-t$, and for $-1 < t \leq 0$
we have $1-t \leq \frac{1}{1+t}$. Hence, the sum of the first two terms is
negative, and $f_t'(x) \leq 0$. By continuity, $f_t$ is decreasing, and so
$q_t(n)$ is decreasing. Since by (4.16) this sequence is also bounded, it is
convergent and passing to limit in (4.15) we get that this limit is $1$ as
claimed.  \qed

Some of our results on orthogonal polynomials, in particular those concerning
Legendre and Laguerre polynomials, have been already known (see, e.g. 
\cite{eweida}), but here we derived them in a simple and unified manner,
almost automatically.

So far we have been considering only the two-term recurrences. The ``calculus''
method, however, works as well for higher order recurrences, as the following
example shows.

Let $L(n)$ be the number of graphs on the vertex set $[n]$, whose every
component is a cycle, and define $L(0):=1$. Then $L(1)=1$, $L(2)=2$, $L(3)=5$,
$L(4)=17$ etc. The following recurrence holds (\cite{stanleyII}, Ex. 5.22):
\be
L(n)=nL(n-1)-\binom{n-1}{2}L(n-3), \quad n\geq 3.
\ee

{\bf Theorem 4.17}\\ The sequence $\left ( L(n) \right ) _{n\geq 0}$ is 
log-convex.

{\bf Proof}\\ Let $q(n) = L(n)/L(n-1)$. From (4.17), dividing with $L(n-1)$,
we obtain 
\be
q(n)=n-\binom{n-1}{2}\frac{1}{q(n-1)q(n-2)}, \quad n\geq 3,
\ee
with the initial conditions $q(1)=1$, $q(2)=2$. The value $q(3)$ is equal to
$5/2$. We claim that $\left ( q(n) \right ) _{n \geq 1}$ is an increasing 
sequence. To this end, define the function $f:[1,\infty) \rightarrow \R$ by
\be
f(x)= \cases {
x &, if $x\in [1,2]$,\cr
\frac{x+2}{2} &, if $x\in [2,3]$,\cr
x- \frac {(x-1)(x-2)}{2}\frac {1}{f(x-1)f(x-2)} &, if $x\geq 3$.\cr
}
\ee
The function $f$ is obviously continuous, rational on any $[n,n+1]$, $n \in \N$,
without poles on $[n,n+1]$ and clearly increasing on $[1,3]$. Let us first 
prove an a priori bound
\be
f(x) \geq \sqrt {x+1}, \quad x \geq 2.
\ee
This is clearly true for $x\in [2,3]$. Suppose that it is true for $x\in [2,n]$,
and let $x\in [n,n+1]$. Then from (4.20) we have
$$f(x) \geq x-\frac{(x-1)(x-2)}{2}\frac{1}{\sqrt {x}\sqrt {x-1}} = 
x-\frac{(x-2)\sqrt{x-1}}{2\sqrt{x}},$$
and it is easy to check that $$x-\frac{(x-2)\sqrt{x-1}}{2\sqrt{x}} \geq
\sqrt {x+1}, \quad x \geq 2.$$ The function $f$ is differentiable for all $x 
\geq 1$, $x\not\in \N$. Let $x\in (n,n+1)$ for some $n \geq 2$ and denote
$f_1=f(x-1)$, $f_2=f(x-2)$. Suppose inductively that $f_1' \geq 0$, $f_2'
\geq 0$. Then
$$f'(x)=1-\frac{2x-3}{2}\frac{1}{f_1f_2}+\frac{(x-1)(x-2)}{2}\frac{f_1'f_2
+f_1f_2'}{(f_1f_2)^2}.$$
The last term is positive, and from (4.21) it follows that
$$1-\frac{2x-3}{2}\frac{1}{f_1f_2} \geq 0 \Longleftrightarrow
f_1f_2 \geq \sqrt {x}\sqrt {x-1}=\sqrt{x^2-x}\geq x-\frac{3}{2},$$
for $x \geq 2$. Hence, $f'(x) \geq 0$ and theorem is proved. \qed

Consider now a general three-term recurrence of the form
\be
a(n)=R(n)a(n-1)+S(n)a(n-2)+T(n)a(n-3), \quad n \geq 4.
\ee
Dividing this recurrence by $a(n-1)$ and denoting $q(n)=a(n)/a(n-1)$, we obtain
the recurrence for $q(n)$:
\be
q(n)=R(n)+\frac{S(n)}{q(n-1)}+\frac{T(n)}{q(n-1)q(n-2)}, \quad n \geq 3,
\ee
with some given initial conditions, say, $q(1)=b_1$, $q(2)=b_2$, $b_2>b_1>0$.

Suppose we want to prove that $(q(n))_{n \geq 1}$ is an increasing sequence.
Again, form the function $f:[1,\infty ) \rightarrow \R$ mimicking the rule
(4.22) and starting with the linear functions on $[1,2]$ and $[2,3]$, 
connecting the points $(1,b_1)$ and $(2,b_2)$ and the points $(2,b_2)$ and 
$(3,b_3)$, respectively, where $b_3=q(3)=R(3)+\frac{S(3)}{b_2}+\frac{T(3)}
{b_1b_2}$ (supposing, of course, $b_3 \geq b_2$). For $x\geq 3$, the function $f$ is defined by replacing $q$ by 
$f$ and $n$ by $x$ in (4.22).
\be
f(x)=R(x)+\frac{S(x)}{f(x-1)}+\frac{T(x)}{f(x-1)f(x-2)},
\ee
for $x \geq 3$ (or $x \geq n_0+1$, for some $n_0 \in \N$). For a fixed $x$ we
write simply $f(x)=f$, $R(x)=R$, $S(x)=S$, $T(x)=T$, $f(x-1)=f_1$, $f(x-2)=
f_2$, etc. So, we write (4.23) simply as
\be
f=R+\frac{S}{f_1}+\frac{T}{f_1f_2} \Longleftrightarrow f_1f_2f=Rf_1f_2+Sf_1+T.
\ee
We assume that $R$, $S$ and $T$ are ``good enough'' functions, in the sense
that $f$ is differentiable on open intervals $(n,n+1)$ for integers $n \geq 2$
(or $n \geq n_0$). This is usually a consequence of an a priori bound of the
type
\be
0 < m(x) \leq f(x) \leq M(x),
\ee
for some well-behaved functions $m(x)$ and $M(x)$.

Taking the derivative $d/dx$ of both sides of (4.24) gives us
\be
f'=\frac{1}{f_1f_2}\left [ F-\frac{Sf_2+T}{f_1}f_1'-\frac{T}{f_2}f_2' \right ] ,
\ee
where
\be
F=R'f_1f_2+S'f_2+T'.
\ee
Suppose that $R, S, T \geq 0$ and $R'$, $S'$, $T' \geq 0$. By substituting the
analogous expressions for $f_1'$ and $f_2'$ in (4.26) and supposing inductively
that $f_2'$, $f_3'$ and $f_4' \geq 0$ we obtain (denoting $F_1=F(x-1)$, $R_1
=R(x-1)$, etc.) that:
$$f_1f_2f'=F-\frac{Sf_2+T}{f_1f_2f_3}F_1-\frac{T}{f_2f_3f_4}F_2+Bf_2'+Cf_3'+
Df_4',$$
where $B$, $C$ and $D$ are non-negative. To prove that $f' \geq 0$, assuming
inductively that $f_2'$, $f_3'$ and $f_4' \geq 0$, it is enough to prove the
following:
$$F-\frac{Sf_2+T}{f_1f_2f_3}F_1-\frac{T}{f_2f_3f_4}F_2 \geq 0,$$
or equivalently
\be
f_1f_2f_3f_4F \geq (Sf_2+T)f_4F_1+Tf_1F_2.
\ee
Taking into account (4.25), it suffices to prove a stronger inequality
\be
m^4F \geq (MS+T)MF_1+MTF_2
\ee
for all $x \geq n_0$, for some $n_0 \in \N$, and to check that $f'(x)\geq 0$
for $x \in (k,k+1)$, $k < n_0$.

Let us now consider a concrete example. A {\bf Baxter permutation} is defined 
in \cite{stanleyII}. The numbers $B(n)$ of Baxter permutations in $\Sigma _n$
satisfy the recurrence
\bea
(n+1)(n+2)(n+3)(3n-2)B(n)&=&2(n+1)(9n^3+3n^2-4n+4)B(n-1) \cr
& &+(3n-1)(n-2)(15n^2-5n-14) B(n-2)\cr
& &+8(3n+1)(n-2)^2(n-3)B(n-3), n\geq 4, 
\eea
together with the initial conditions $B(0)=1$, $B(1)=1$, $B(2)=2$, $B(3)=6$.

{\bf Theorem 4.18}\\ The numbers $B(n)$ of Baxter permutations are log-convex.
The limit $\lim _{n \rightarrow \infty} B(n)/B(n-1)$ exists and is equal to $8$.

{\bf Proof}\\ Let $q(n)=B(n)/B(n-1)$. From (4.30) we form the recurrence for
$q(n)$'s, and then according to the initial values $q(1)=1$, $q(2)=2$, $q(3)=
3$, we form the function $f:[1,\infty) \rightarrow \R$ defined by $f(x)=x$ on
$[1,3]$ and for $x \geq 3$ by the rule
\bea
(x+1)(x+2)(x+3)(3x-2)f(x)&=&2(x+1)(9x^3+3x^2-4x+4)\nonumber \\
& &+\frac{(3x-1)(x-2)(15x^2-5x-14)} {f(x-1)}\\
& &+\frac{8(3x-1)(x-2)^2(x-3)}{f(x-1)f(x-2)}. \nonumber
\eea
This, written in the form of (4.23), yields to conclude that $R,S,T$ are 
positive rational functions with no poles on $[1,\infty )$. We also see that
$R(x) \nearrow 6$, $S(x) \nearrow 15$, $T(x) \nearrow 8$ as $x \rightarrow 
\infty$.

It is apparent from (4.31) that $f(x)$ is a piecewise rational function, i.e.
rational on intervals $[n,n+1]$, $n \in \N$. For example, for $x\in [3,4]$,
$$f(x)=\frac{18x^5+51x^4-122x^3-87x^2+200x+12}{(x-1)(x+1)(x+2)(x+3)(3x-2)}.$$
Clearly, $f(x)$ is continuous everywhere. It can be checked, by induction on
$n$, that the function $f$ is bounded. More precisely, 
\be
7 \leq f(x) \leq 9, \quad x \geq 47.
\ee
We want to prove that $f$ is an increasing function. From the above a priori
bound, it follows that $f$ is differentiable on all open intervals $(n,n+1)$,
$n \in \N$. The non-negativity of $f'(x)$ can be checked for $x \leq 49$,
using e.g. {\it Mathematica}, and for $x \geq 49$, it follows from the stronger
inequality,
$$7^4 F \geq (9 \cdot 15+8)\cdot 9 F_1+ 9\cdot 8 F_2,$$
obtained by substituting appropriate values for $m$, $M$ in (4.29). As this
inequality is true for all $x \geq 9$, the first claim follows. The second
claim follows passing to the limit in (4.30). \qed

Many other three-or-higher-term recurrences can be investigated by this method.
Let us only mention that the four-term recurrence (3.15) for the numbers
$S^{(1)}(n)$ of secondary structures of rank $1$ can also be shown to be
log-convex by this ``calculus'' method. The details (rather tedious) are given
in \cite{doslicphd}, \cite{doslvelj}. Another example is the number $S_n$ of
$n \times n$ symmetric matrices with entries $0$, $1$, $2$, whose sums of all
rows and all columns are equal to $2$. These numbers satisfy (\cite{stanleyII})
the following recurrence
$$S_n=(2n-1)S_{n-1}-(n-1)(n-2)S_{n-2}-(n-1)(n-2)S_{n-3}+\frac{1}{2}(n-1)(n-2)
(n-3)S_{n-4},$$
starting with $S_0=1$, $S_1=1$, $S_2=3$, $S_3=11$. The sequence $\left ( S_n
\right ) _{n \geq 0}$ is also log-convex.

Finally, as it should be clear by now, this method applies to any $P$-recursive
sequence $(a(n))$, satisfying a recurrence of the form
\be
Q(n)a(n)=P_d(n)a(n-1)+P_{d-1}(n)a(n-2)+ \ldots + P_0(n)a(n-d-1), \quad n \geq
d+1,
\ee
where $d \geq 0$ is an integer and $P_0, \ldots ,P_d, Q$ real polynomials, $Q
>0$. The corresponding function for the successive quotients $q(n)=\frac{a(n)}{a(n-1)}$ is given by the functional equation
\be
Q(x)f(x)=P_d(x)+\frac{P_{d-1}(x)}{f(x-1)}+ \ldots + \frac{P_0(x)}{f(x-1) \ldots
f(x-d)}, \quad x\geq d+1.
\ee
As we have seen, the most important thing in this approach is to express the
derivative $f'(x)$ in terms of previous derivatives. So, fix a point $x\in
(n,n+1)$, $n > d$, and write as before for short $f=f(x)$, $f_j=f(x-j)$, $j=1,2,
\ldots $, $P_i=P_i(x)$, $Q_i=Q_i(x)$, $i=0,1, \ldots ,d$. Then (4.34) can be 
written as $$Qf=\sum _{i=0}^d \frac{P_i}{f_1f_2 \ldots f_{d-i}},$$ or by 
denoting the product of all values by $\Pi $, i.e. $\Pi = f_1f_2 \ldots
f_d$, and by $\Pi _j = f_1f_2 \ldots f_j$ the partial products (so $\Pi _d
=\Pi$ and $\Pi _0 =1$), as
\be
\Pi Qf = \sum _{i=o}^d(f_{d-i+1} \ldots f_d)P_i.
\ee
Taking the derivative $d/dx$ of both sides in (4.35), after some manipulations,
we obtain the following formula (in terms of Wronskians):
\be
f'=\frac{1}{Q^2}\sum _{i=0}^d \left | {Q \atop Q'} \quad {{P_i} \atop {P_i'}}
\right | \frac{1}{\Pi _{d-i}}-\frac{1}{Q}\sum _{i=0}^d \frac{P_i \Pi _{d-i}'}{\Pi _{d-i}^2}.
\ee
In particular, for $d=1$ this reduces to
\be
f'=\frac{1}{Q^2f_1}\left | {Q \atop Q'} \quad {{P_0} \atop {P_0'}} \right |+
\frac{1}{Q^2}\left | {Q \atop Q'} \quad {{P_1} \atop {P_1'}} \right |-
\frac{P_0}{Qf_1^2}f_1'.
\ee
With a priori bounds $0<m(x) \leq f(x) \leq M(x)$ and by substituting $f_1'$
in (4.37), one gets almost instant proofs of log-behavior. For example, if $P_0
\leq 0$ and if we want to prove the log-convexity, hence assuming $f_1' \geq 0$, then if the first Wronskian $W_0$ in (4.37) is positive, we only have to check
$\frac{1}{M}W_0+W_1 \geq 0$, and check that $f$ increases at the beginning.


Of course, not every (combinatorially relevant) sequence satisfying a recurrence
of this type can be expected to have a reasonable log-behavior; it is enough to
recall here the sequences $e_k(n)$ from Section 2, whose log-behavior is rather
chaotic for $k \geq 3$.

Let us only mention here that the log-convexity of secondary
structure numbers of general rank $l$ can also be proved by calculus method,
using the explicit formulae from Proposition 3.9 and formula (4.36). The details
will appear elsewhere.

As a finall remark, note that our approach applies also to linear 
nonhomogeneous recurrences for positive numbers. So, for example, let $(a(n))$
be given by the linear recurrence of the first order
\be
a(n)=R(n)a(n-1)+S(n).
\ee
Consider the quotients $q(n)=\frac{a(n)}{a(n-1)}$ and note that
\be
a(n)=q(n)q(n-1) \ldots q(2)a(1), \quad n\geq 2.
\ee
Then, dividing (4.38) by $a(n-1)$ we obtain a (long) recurrence for $q(n)$'s:
\be
q(n)=R(n)+\frac{S(n)}{q(n-1)q(n-2) \ldots q(2)a(1)}.
\ee
To get a short recurrence for $q(n)$'s, substitute for $a(n)$ and $a(n-1)$ the
corresponding products (4.39) in (4.38) ($n\geq 3$):
$$\displaylines{
\qquad q(n)q(n-1)\ldots q(2)a(1)=R(n)q(n-1)\ldots q(2)a(1)+S(n)= \hfill \cr
\hfill R(n)\frac{q(n)q(n-1)\ldots q(2)a(1)}{q(n)}+S(n)\frac{q(n)q(n-1)\ldots q(2)a(1)}{q(n)q(n-1)\ldots q(2)a(1)}. \qquad\cr
}$$
From there we get
$$\frac{1}{q(n)q(n-1)\ldots q(2)a(1)}=\frac{1}{S(n)}\left [ 1-\frac{R(n)}{q(n)} \right ],$$
and then
\be
\frac{1}{q(n-1)\ldots q(2)a(1)}=\frac{1}{S(n-1)}\left [ 1-\frac{R(n-1)}{q(n-1)} \right ] .
\ee
Substituting (4.41) in (4.40) yields a short recursion for $q(n)$'s:
$$q(n)=R(n)+\frac{S(n)}{S(n-1)}-\frac{R(n-1)S(n)}{S(n-1)}\frac{1}{q(n-1)}.$$
Similarly, for a second order linear recurrence
$$a(n)=R(n)a(n-1)+S(n)a(n-2)+T(n),$$
we obtain
$$q(n)=R(n)+\frac{S(n)}{q(n-1)}+\frac{T(n)}{T(n-1)}\left [ 1- \frac{R(n-1)}{q(n-1)}-\frac{S(n-1)}{q(n-1)q(n-2)} \right ]. $$
Then we can proceed as before.

\section{Calculus method in two variables}

\setcounter{equation}{0}
We shall outline our method for non-negative sequences $a(n,k)$ in two integer
variables $n,k \geq 0$ (or $n \geq n_0$, $k \geq k_0$). Suppose (as often in
combinatorics) that the numbers $a(n,k)$ satisfy a two-term recurrence
of the form
\be
a(n,k)=R(n,k)a(n-1,k-1)+S(n,k)a(n-1,k),
\ee
with some known functions $R$ and $S$, together with some initial values, usually
of the type $a(0,0)=a$, $a(1,0)=b$, $a(1,1)=c$. Suppose we want to prove that
the sequence $a(n,k)$ is log-concave in $k$, i.e. that $a(n,k)^2 \geq a(n,k-1)
a(n,k+1)$, for all $n,k$. Here is what we do. Write down (5.1) with $k$ 
replaced by $k-1$:
\be
a(n,k-1)=R(n,k-1)a(n-1,k-2)+S(n,k-1)a(n-1,k-1).
\ee
Denote
\be
q(n,k)=\frac{a(n,k)}{a(n,k-1)},
\ee
and divide (5.1) by (5.2) (always assuming we do not divide by zero).
\begin{eqnarray*}
q(n,k)&=&\frac{R(n,k)a(n-1,k-1)+S(n,k)a(n-1,k)}{R(n,k-1)a(n-1,k-2)+S(n,k-1)a(n-1,k-1)}
=\frac{R(n,k)+S(n,k)q(n-1,k)}{\frac{R(n,k-1)}{q(n-1,k-1)}+S(n,k-1)}\\
&=& q(n-1,k-1)\frac{R(n,k)+S(n,k)q(n-1,k)}{R(n,k-1)+S(n,k-1)q(n-1,k-1)}.
\end{eqnarray*}
Equivalently,
\be
q(n,k)[R(n,k-1)+S(n,k-1)q(n-1,k-1)]=q(n-1,k-1)[R(n,k)+S(n,k)q(n-1,k)].
\ee
The log-concavity of $a(n,k)$'s is equivalent to $q(n,k) \geq q(n,k+1)$, for
any fixed $n$ and all $k$. The idea is again to pass to a ``continuation'' of
(5.4) by letting $n \rightarrow x$, $k \rightarrow y$, $q \rightarrow f$ and
obtaining the functional equation
\be
f(x,y)[R(x,y-1)+S(x,y-1)f(x-1,y-1)]=f(x-1,y-1)[R(x,y)+S(x,y)f(x-1,y)].
\ee
We assume that $R$ and $S$ are ``good enough'' functions, in the sense that $f$
is continuous everywhere and smooth on open cells $(n,n+1) \times (m,m+1)$, for all $m,n$.
What we want to prove is that $f$ is decreasing in $y$ for any fixed $x$. Fix a
point $(x,y)$ in an open cell $Q=(n,n+1) \times (m,m+1)$, and prove inductively
that \be
\frac{\partial f}{\partial y}(x,y) \leq 0.
\ee
For the fixed pair $(x,y)$ write for short $f_{ij}=f(x-i,y-j)$, for $i,j=0,1,
2,\ldots $, and similarly for $R$ and $S$. So, $f_{00}=f(x,y)$, $R_{01}=
R(x,y-1)$ etc. In this notation, (5.5) can be written as
\be
f[R_{01}+S_{01}f_{11}]=f_{11}[R+Sf_{10}].
\ee
Now take the partial derivative $\partial /\partial y$ of both sides in (5.7).
We have
$$\frac{\partial f}{\partial y}\left [ R_{01}+S_{01}f_{11} \right ] +
f\left [ \frac{\partial R_{10}}{\partial y}+f_{11}\frac{\partial S_{01}}
{\partial y}+S_{01}\frac{\partial f_{11}}{\partial y} \right ] =
\frac{\partial f_{11}}{\partial y}\left [ R+Sf_{10}\right ]+f_{11}
\left [ \frac{\partial R}{\partial y}+f_{10}\frac{\partial S}{\partial y}
+S\frac{\partial f_{10}}{\partial y}\right ].$$
$$\frac{\partial f}{\partial y}\left [ R_{01}+S_{01}f_{11} \right ]=
f_{11}\left [ \frac{\partial R}{\partial y}+f_{10}\frac{\partial S}{\partial y}
\right ] - \left [ \frac{\partial R_{10}}{\partial y}+f_{11}\frac{\partial S_{01}}{\partial y}\right ]f
+\frac{\partial f_{10}}{\partial y}f_{11}S+\frac{\partial f_{11}}{\partial y}
\left [ R+Sf_{10}-S_{01}f\right ].$$
Substituting here $f$ from (5.7), we get
\bea
\frac{\partial f}{\partial y}\left [ R_{01}+S_{01}f_{11} \right ]&=& 
\frac{f_{11}}{R_{01}+S_{01}f_{11}}\left [ (R_{01}+S_{01}f_{11})(\frac{\partial R}{\partial y} +f_{10}\frac{\partial S}{\partial y})-
(R+Sf_{10})(\frac{\partial R_{01}}{\partial y}+f_{11}\frac{\partial S_{01}}{\partial y} \right ] \cr
& &+f_{11}S\frac{\partial f_{10}}{\partial y}+\frac{(R+Sf_{1,0})R_{01}}{R_{01}+S_{01}f_{11}}\frac{\partial f_{11}}{\partial y}.
\eea
Assume that $R$ and $S$ are positive. Hence if $f$ starts with some positive
values, then $f$ can be considered positive, too. Suppose inductively that 
$$ \frac{\partial f_{10}}{\partial y}, \frac{\partial f_{11}}{\partial y}
\leq 0.$$
Then the last two terms in (5.8) are negative and to prove (5.6), it is 
enough to prove that the first term is negative, too. In other words, to
conclude (inductively) that (5.6) holds, it is enough to prove that the ``free''
term is non-positive, i.e. that
\be
F=(R_{01}+S_{01}f_{11})\left ( \frac{\partial R}{\partial y}+f_{10}\frac{\partial S}{\partial y}\right ) - (R+Sf_{10})\left (\frac{\partial R_{01}}{\partial y}+f_{11}\frac{\partial S_{01}}{\partial y} \right ) \leq 0.
\ee
So, if we can check that $f$ begins decreasingly in $y$ and assuming inductively
that $f$ is decreasing in $y$, then by (5.9) we can conclude that $f$ is 
decreasing in $y$ at the point $(x,y)$ and then, by continuity, that $f$ is
decreasing in $y$ everywhere.

Note that the following inequalities imply (5.9) (simply by comparing similar
terms):
$$R_{01}\frac{\partial R}{\partial y}\leq R\frac{\partial R_{01}}{\partial y},
S_{01}\frac{\partial R}{\partial y}\leq R\frac{\partial S_{01}}{\partial y},
R_{01}\frac{\partial S}{\partial y}\leq S\frac{\partial R_{01}}{\partial y},
S_{01}\frac{\partial S}{\partial y}\leq S\frac{\partial S_{01}}{\partial y}.$$
In terms of Wronskians, writing $G'$ for $\frac{\partial G}{\partial y}$,
these inequalities can be written in the form:
\be
\left | {R_{01} \atop R_{01}'} \quad {R \atop R'} \right | \leq 0,
\left | {S_{01} \atop S_{01}'} \quad {R \atop R'} \right | \leq 0,
\left | {R_{01} \atop R_{01}'} \quad {S \atop S'} \right | \leq 0,
\left | {S_{01} \atop S_{01}'} \quad {S \atop S'} \right | \leq 0.
\ee
Instead of formalizing everything (which can be done with a little care), let 
us take an example.

{\bf Example 5.1} \\
The {\bf Eulerian number} $E(n,k)$, is the number of permutations $\pi $ from
$\Sigma _n$ with exactly $k$ ascents, i.e. with exactly $k$ places where
$\pi _j < \pi _{j+1}$.
We know (see, e.g. \cite{gkp}) that these numbers satisfy the recurrence
$$E(n,k)=(n-k)E(n-1,k-1)+(k+1)E(n-1,k),$$
with the initial conditions $E(0,k)=\delta _{0k}$, $E(n,0)=1$, $n,k \geq 0$.
Then (5.4) becomes
$$q(n,k)[n-k+1+kq(n-1,k-1)]=q(n-1,k-1)[n-k+(k+1)q(n-1,k)]$$
with $q(n,k)=E(n,k)/E(n,k-1)$,
for $1\leq k \leq n$. The initial conditions are $q(2,1)=1$, $q(3,1)=4$,
$q(3,2)=1/4$ and (we extra define) $q(2,k)=0$, for $k \geq 2$. We want to prove
that $q(n,k) \geq q(n,k+1)$, for any fixed $n$ and all $k$. In the sense of the
above discussion and notations, here we have $R(n,k)=n-k$, $S(n,k)=k+1$. Passing to
the natural ``continuation'', (5.7) becomes
\be
f[x-y+1+yf_{11}]=f_{11}[x-y+(y+1)f_{10}].
\ee
In fact, we define the function $f :[2,\infty ) \times [1,\infty )\rightarrow 
\R$ first on two shaded strips in Fig. 6 below and then continue by the rule
(5.11). In Fig. 6 we indicated the values $f(n,k)=q(n,k)$ for $n \geq 2$, 
$k \geq 1$ in
the lattice nodes. On the vertical walls of the (square) cells $Q_2$, $Q_3$,
$Q_4, \ldots $, as well as on their lower horizontal walls, we define $f$ to be
appropriate linear functions.  The upper walls of $Q_3$, $Q_4, \ldots $, are
determined by (5.11). The right walls of $Q_3'$, $Q_4', \ldots $ are also 
determined by (5.11). We fill in $f$ on cells $Q_n$, $Q_n'$ by appropriate
homotopies connecting (possibly nonlinear) functions on the walls. 
\begin{figure}[h] \centerline {
\epsfig{file=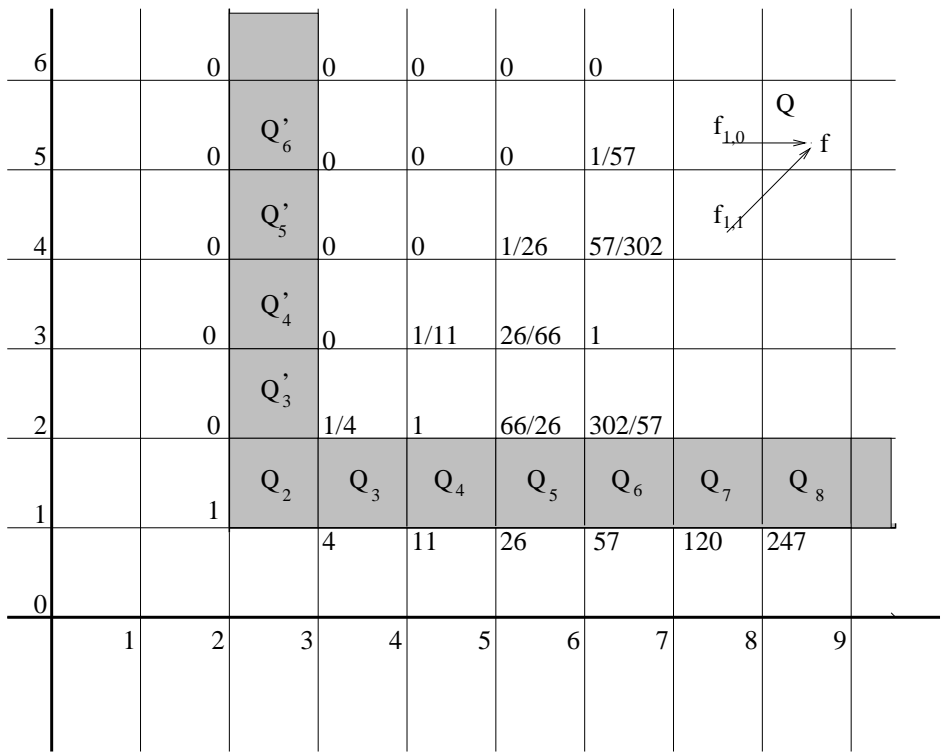,height=10cm,width=12.5cm,silent=}}
\caption{The boundary values of function $f$}
\label{sec}
\end{figure}

So, let $f(x,2)=2^n(x-n)+2^n-n-1$, for $x\in [n,n+1]$, $n\geq 2$, 
$f(2,y)=2-y$ for $y \in [1,2]$ and $f(2,y)=0$ for $y \geq 2$. Further, let
$f(3,y)=\frac{31-15y}{4}$ for $y \in [1,2]$, and, in general, let $f(n,y)$, 
$y \in [1,2]$, be the linear function in $y$ between the points $(n,q(n,1))$
and $(n,q(n,2))$. The cell $Q_2$ is surrounded by linear functions
$a(x)=3x-5$, $x \in [2,3]$, $b(y)=\frac{31-15y}{4}$, $y \in [1,2]$, $c(x)=
\frac{x-2}{4}$, $x \in [2,3]$ and $d(y)=2-y$, $y \in [1,2]$. Extend $f$ to 
$Q_2$ by the homotopy
$$f(x,y)=(2-y)a(x)+(y-1)c(x)=\frac{1}{4}(23x+18y-11xy-38 ), \quad (x,y)\in Q_2.$$
The cell $Q_3'$ is surrounded by functions $a_3'(x)=\frac{x-2}{4}$, 
$x \in [2,3]$, $b_3'(y)=\frac{(3-y)^2}{4-y+y(3-y)}$, $y\in[2,3]$, $c_3'(x)=0$, 
$x \in [2,3]$ and $d_3'(y)=0$, $y\in[2,3]$. Extend $f$ to $Q_3'$ by the
homotopy
$$f(x,y)=(3-x)d_3'(y)+(x-2)b_3'(y)=\frac{(x-2)(3-y)^2}{4-y+y(3-y)}, \quad (x,y)
\in Q_3'.$$
Extend $f$ to $Q_4'$, $Q_5', \ldots $ to be zero. Next, extend $f$ to $Q_3$, $Q_4, \ldots$
by the homotopies connecting the lower walls of $Q_n$, given by $a_n(x)=
2^n(x-n)+2^n-n-1$, for $x\in [n,n+1]$, $n\geq 2$, and the upper walls,
given by rational functions $c_n(x)$ determined  inductively on $n$ by (5.11),
thus obtaining $f \Big | _{Q_n}$ by
$$f(x,y)=(2-y)a_n(x)+(y-1)c_n(x), \quad (x,y) \in Q_n.$$
For example, since $a_2(x)=3x-5$, $a_3(x)=7x-17$, $c_2(x)=\frac{x-2}{4}$, then
$$c_3(x)=\frac{a_2(x-1)[x-y+(y+1)c_2(x-1)]}{x-y+1+ya_2(x-1)}=\frac{3x-8}{4},$$
and hence $f \Big | _{Q_3}$ is given by
$$f(x,y)=(2-y)a_3(x)+(y-1)c_3(x)=\frac{53x+60y-25xy-128}{4}.$$
In this way, $f$ is well defined on the shaded strips on Fig. 6 and extended
to $[2,\infty )\times [1,\infty )$ by the rule (5.11). 

It is easy to check that $f$ is continuous and nonnegative and that it is a 
rational function with no poles on any open cell, and hence smooth on any open
cell. It is also easy to check inductively on $n$ that $f$ is decreasing in
$y$ on $Q_n$ and $Q_n'$ for $n \geq 2$, i.e. $\frac{\partial f}{\partial y}
(x,y) \leq 0$ for $(x,y) \in int(Q_n) \cup int(Q_n')$.

Now that we have elaborated carefully the ``beginning'' of $f$, the rest is
more-or-less automatic. The inequality (5.9) reduces to $(x-y+1+yf_{11})(-1+
f_{10}) \leq (x-y+yf_{10}+f_{10})(-1+f_{11})$. This is equivalent to
$$-1+xf_{10}+f_{10} \leq -f_{10}+xf_{11}+f_{10}f_{11},$$ or, after some 
rearrangement, $$f_{11}f_{10}-2 f_{10}+1+x(f_{11}-f_{10}) \geq 0. $$
The second term is non-negative by the induction hypothesis, and the rest is
non-negative since $$f_{11}f_{10}-2 f_{10}+1 \geq f_{10}^2-2 f_{10}+1 =
(f_{10}-1)^2 \geq 0.$$
Hence, $\frac{\partial f}{\partial y}(x,y) \leq 0$ for $(x,y) \in int(Q)$,
for all $Q$. So, $f$ is decreasing in $y$ on every open cell, and hence by
continuity, $f$ is decreasing in $y$ everywhere. In particular, $q(n,k) =
f(n,k) \geq f(n,k+1) = q(n,k+1)$ and we are done.

In the same manner we can prove that $\frac{\partial f}{\partial x}(x-y+1)
\geq f$, hence $\frac{\partial f}{\partial x} \geq 0$, for any $(x,y) \in int(Q)$, for any $Q$, and this implies $q(n+1,k) \geq q(n,k)$, for any fixed $k$
and all $n$. Hence,
$$\frac{E(n+1,k)}{E(n+1,k-1)} \geq \frac{E(n,k)}{E(n,k-1)} \Longleftrightarrow
\left | {{E(n+1,k)} \atop {E(n,k)}} \quad {{E(n+1,k-1)} \atop {E(n,k-1)}}
\right | \geq 0.$$

Apart from settling the ``beginning'' of $f$, we can (almost automatically now)
prove the well-known log-concave behavior in the second variable of the
binomial coefficients (and in general find the log-concave behavior when $R$ 
and $S$ in (5.1) are constants), $q$-binomial coefficients, Stirling numbers
of the first and second kind, Schl\"afli numbers, cover many particular results (e.g 
\cite{kurtz}) and so on.

Of course, the method can be extended in a few ways; for example, to 
three-or-more term recurrences, recurrences for three or more variables
$a(n,k,l)$ etc. but we shall not consider it here.

In a word, a general idea of this method is as follows. Combinatorics gives a 
recurrence. Pass to the quotients of the neighboring members (in the variable
under consideration), pass to the natural ``continuation'' $f$, find some
bounds (upper, lower or both, depending on the nature of the problem) of $f$,
make sure that $f$ is differentiable on open cells, examine the rate of growth
of $f$ at the ``beginning'' (i.e. check the sign of the derivative there), and
finally, prove inductively from the associated functional equation corresponding
to the recurrence, that the sign of the derivative remains the same.

\end{document}